\documentstyle[11pt,fleqn]{article}

\begin{document}
\setcounter{page}{1}
\thispagestyle{empty}

\vspace*{10mm}
\noindent
\centerline{{\Large  Linear representations of groups  with }}
\vskip 2mm
\centerline{{\Large  translation  invariant  defining relationships. }}
\vskip 2mm
\centerline{{\Large  Some new series of braid  group representations }}
\vskip 2mm
\centerline{{\Large  and new invariants of links and knots }}

\vskip 6mm
\noindent
\centerline{Vladimir K. Medvedev (Kiev University)}
\vskip 10mm

In this paper we indicate one method of construction of linear
representations of groups and algebras with translation invariant
(except, maybe , finite number) defining relationships. As an
illustration of this method, we give one approach to the construction of linear
representations of braid group and derive some series of such
representations.
Some invariants of oriented knots and links are constructed.
The author is grateful to Yuri Drozd, Sergey Ovsienko and other members
 of algebraic seminar at Kiev University for the
creative atmosphere without which this work could hardly appear.

\vskip 6mm
\noindent
{\Large\bf 1} \ {\large {\bf The description of the approach on the example
of the
\hspace *{6mm} braid group}}
\vskip 3mm

\noindent 1.1. {\bf Definition.} \ {\it Let $N$ be the set of the natural
numbers.
Braid group $B_{\infty} (B_m)$ is the group with generators $t_i$,
$i\in N \, \ (i\in \{1,...,m-1\})$ and
\linebreak the following defining relations:  }
\begin{displaymath}
\matrix {
({\rm i}) \ t_i t_{i+1} t_i =t_{i+1} t_it_{i+1},\hfill\cr\cr
({\rm ii}) \ t_it_j=t_jt_i \ when  \ |j-i|> 1.\hfill }
\end{displaymath}
\vskip 2mm

\noindent
1.2. \ For natural $n,\lambda, \lambda < n$ let $M_1 =\{1,2,...,n\} \subset N$
and
for any $i\in N $ let the set $M_{i+1}$ to be obtained of the set $M_i$ by the
right translation on the number $n-\lambda$.
It is obvious
that $\vert {M_i}\vert =n$ and $\vert {M_i \cap M_{i+1}}\vert =\lambda$.

We shall call linear operator $f: V\to V$ on the vector space V over the field
$k$
nondegenerate or invertible if
there exists linear operator $g : V\to V$ such that $fg =gf =I_V$.
For the family $\{V_s\}_{s\in N}$ of the vector spaces over the field $k$ let
$W ={\mathop {\oplus } \limits_{s\in N}} V_s$ and
$V_M ={\mathop {\oplus } \limits_{s\in M}} V_s$ for any $M\subset N$.
Let  for any $i\in N$ there is some nondegenerate linear operator
$ T_i$ on the vector space $V_{M_i}$. Operator $ T_i$ is defined by the
matrix $\{T_i  (s_1 +1\ $-- $\min M_i$, $s_2 +1\ $-- $\min M_i) \}
_{\scriptstyle s_1 \in M_i \atop\scriptstyle s_2 \in M_i} $, where
$T_i (s_1+1\ $-- $\min M_i, s_2+1\ $-- $\min M_i)$ is the linear operator from
the
vector space $V_{s_2}$ into the vector space $V_{s_1}$.
We shall extend this operator to the linear
operator $\tilde t_i$ on $W$ by setting
$\tilde t_i=(T_i)_{V_{M_i}} \oplus1_{V_{N\setminus M_i}}$.
Relations \mbox{1.1. (i)--(ii)}
that hold  if and only if they  hold on the vector spaces $V_{M_i\cup M_{i+1}}$
(for 1.1. (i))
and $V_{M_i\cup M_j}$ (for 1.1. (ii)), define some relations for the
matrix elements
of matrices $T_i$. If they are true we get linear representation of the
braid group $B_{\infty}$.
We note that if $\lambda \le {n\over 2}$ relations
1.1. (ii) hold automatically.
\vskip 3mm

{\bf Remark.} \ {\it If there exists $d\in N$ such that for any $i$
$\ T_i =T_{i+d}\ $ in
some basises in $V_s, s\in N$ (that is if we consider periodic case with
period $d$), then we have already the finite number of
matrices $T_1,...,T_d$ with finite number of relations for their matrix
elements
(and thus, some "algebra with finite number of generators and defining
relations")
such that if they hold we obtain a linear
representation of the group  ${\cal B}_{\infty}$.
The same
approach is true to any group (and some rings) that can be defined by
translation invariant
relationships. It is the crux   of this article.    }
\vskip 6mm
%
  1
The same approach is true if we consider instead of direct sum tensor
product. Namely, for the braid group $B_m$ let us consider
$$
W= \otimes {\mathop {\prod} \limits_{s=1}^{s=\lambda +(m-1) (n-\lambda)}} V_s,
\qquad V_M=\otimes {\mathop {\prod} \limits_{s\in M}} V_s,
$$
where $M\subset N$.
Let $T(i)$, for any $i\in \{1,...,m-1\}$
be some nondegenerate linear operator on the vector space
$V_{M_i}$. Thus $T(i)$ is defined by tensor
$T(i) \in V_{M_i} \otimes V_{M_i}^*$.
We shall define the operator $\tilde t_i$ on vector space $W$ by setting
$\tilde t_i =T(i)_{V_{M_i}} \otimes 1_{V_{ \{1,...,\lambda +(m-1)(n-\lambda)\}
\setminus M_i}}$.
\vskip 5mm

\noindent
Relations 1.1. \  (i)-(ii) that hold if and only if they hold on the vector
space $V_{M_i \cup M_{i+1}}$ (for 1.1 (i)) and $V_{M_i \cup V_{M_j}}$
(for 1.1 (ii))
define some relations for the
$T(i)$. If $\lambda \le {n\over 2}$ relations 1.1 (ii) hold automatically.
The obvious analog of the remark above is true for the case of tensor
products (that is we separate periodic case, when for some
$d\in N$ $\ T(i+d)=T(i)\ $
for any $i$ and, thus, all relations for $T(i)$ are consequences of
the finite number of relations for $T(1), ..., T(d)$ ).
\vskip 6mm

\noindent
{\Large\bf 2} \ {\large  {\bf Linear representations of the braid group}}
\vskip 4mm

\noindent
To illusrate as told above we consider at first the simplest case
\linebreak $(n,\lambda )=(2,1)$.
\vskip 4mm

\noindent
2.1. \ Relations 1.1. (ii) hold  evidently. On the vector space
$V_{M_i \cup M_{i+1} } =V_i \oplus V_{i+1} \oplus V_{i+2}$ matrices
of the linear operators $\tilde t_i, \tilde t_{i+1}$ have the following form:
\begin{displaymath}
T_i =\left( \matrix{ A_i & B_i & 0\cr
C_i & D_i & 0\cr
0 & 0 &1 } \right),\qquad
T_{i+1} =\left( \matrix{1 & 0 & 0\cr
0 & A_{i+1} & B_{i+1} \cr
0 & C_{i+1} & D_{i+1} } \right).
\end{displaymath}

The  condition $\tilde t_i \tilde t_{i+1} \tilde t_i =\tilde t_{i+1} \tilde t_i
\tilde t_{i+1}$,
$i\in N$ holds then and only then when $T_i T_{i+1} T_i =T_{i+1} T_i T_{i+1}$.
Simple calculation shows that it is equivalent to the following relations:
\begin{displaymath}
\matrix{ ({\rm i}) &  A^2_i +B_i A_{i+1} C_i =A_i,\hfill\cr
({\rm ii}) & D^2_{i+1} +C_{i+1} D_i B_{i+1} =D_{i+1},\hfill\cr
({\rm iii}) &  A_{i+1} C_i -C_i A_i =D_i A_{i+1} C_i,\hfill\cr
({\rm iv}) & B_i A_{i+1} -A_i B_i =B_i A_{i+1} D_i ,\hfill\cr
({\rm v}) & C_i B_i -B_{i+1} C_{i+1} =A_{i+1}D_i A_{i+1} -D_i A_{i+1}
D_i,\hfill\cr
({\rm vi}) & C_{i+1} D_i -D_{i+1} C_{i+1} =C_{i+1} D_i A_{i+1},\hfill\cr
({\rm vii}) & D_i B_{i+1} -B_{i+1} D_{i+1} =A_{i+1} D_i B_{i+1}\hfill\cr  }
\end{displaymath}
%
  2  

Let us consider the case of tensor product. Relations 1.1 (ii) hold
evidently. Simple calculation shows that relations
$\tilde t_i \tilde t_{i+1} \tilde t_i =\tilde t_{i+1} \tilde t_i \tilde
t_{i+1}$
are equivalent to the following relations for the tensors
$T(i) \in (V_i \otimes V_{i+1}) \otimes (V_i \otimes V_{i+1} )^*$,
$i=1,...,m-1$:
\begin{displaymath}
({\rm viii}) {\mathop {\sum} \limits_ {k_1, k_2, k_3}}
T(i)^{i_1 i_2}{}_{k_1k_3} \times T(i+1) ^{k_3 i_3}{}_{k_2j_3}
\times T(i) ^{k_1 k_2}{}_{j_1j_2} =
\end{displaymath}
\begin{displaymath}
={\mathop {\sum} \limits_  {k_1,k_2,k_3}} T(i+1)^{i_2 i_3}{}_{k_1k_3} \times
T(i)^{i_1 k_1}{}_{j_1k_2} \times T(i+1)^{k_2k_3}{}_{j_2j_3}
\end{displaymath}
for any $i_1, i_2, i_3, j_1, j_2, j_3, \ 1\le i_1(j_1) \le n_i$,
$1\le i_2 (j_2) \le n_{i+1}$, $1\le i_3 (j_3) \le n_{i+2}$,
where $n_i ={\rm dim}\, V_i$. Here elements $T(i)^{uv}{}_{pq} $
of the field $k$ are defined by some choice of the basis
$\{ e^{(i)}_s \}_{s=1,...,n_i}$ in every $V_i$ and by the equality
\begin{displaymath}
T(i) =\sum T(i) ^{uv}{}_{pq} e^{(i)}_u \otimes e_v^{(i+1)} \otimes e^{(i)p *}
\otimes e^{(i+1)q*}.
\end{displaymath}


Thus we have the following theorem.
\vskip 6mm

{\bf Theorem.} \ {\it Let linear operator $ T_i$ on the vector space
$V_i \oplus V_{i+1}$ with matrix $\left( \matrix{ A_i & B_i\cr C_i &
D_i\cr}\right)$
is nondegenerate for any $i\in N$ and relations
\linebreak {\rm (i)-(vii)} hold.
Then linear operators $T_i$ define representation of the braid
group $B_{\infty}$.

Let $T(i) \in (V_i \otimes V_{i+1} ) \otimes (V_i \otimes V_{i+1} )^*$
be such linear operator on the vector space $V_i \otimes V_{i+1}$ that
matrix $(T(i)^{uv}{}_{pq})$is nondegenerate
for any $i \in \{1,..., m-1\}$ and relations {\rm (viii)} hold.
Then linear operators $ T_i$
define representation of the braid group $B_{m}$.  }
\vskip 4mm

Now let us come to the periodic case according to the approach pointed
out in remark 1.2.

\noindent
\vskip 8mm

\noindent
2.2. \ Let  ${\cal B}$ be an algebra over field $k$ with generators $A, B, C,
D$ and
the following defining relations:
\begin{displaymath}
\matrix{
({\rm i}) &  A^2 +B A C =A,\hfill\cr
({\rm ii}) & D^2 +C D B =D,\hfill\cr
({\rm iii}) &  CB -BC=ADA-DAD,\hfill\cr
({\rm iv}) & B A -A B =B A D ,\hfill\cr
({\rm v}) & AC -CA =DAC,\hfill\cr
({\rm vi}) &  DB -BD =ADB,\hfill\cr
({\rm vii}) & CD -DC=CDA.\hfill\cr  }
\end{displaymath}
%
 3
Let $T\in V^{\otimes 2 } \otimes (V^*)^{\otimes 2}$ be the tensor of type
(2.2),
such that the following relation holds:
\begin{displaymath}
({\rm viii})
 {\mathop {\sum} \limits_ {k_1, k_2, k_3}}
T^{i_1 i_2}{}_{k_1 k_3} \times T ^{k_3 i_3}{}_{k_2 j_3}
\times T ^{k_1 k_2}{}_{j_1j_2} =
\end{displaymath}
\begin{displaymath}
={\mathop {\sum} \limits_  {k_1,k_2,k_3} } T^{i_2 i_3}{}_{k_1 k_3} \times
T^{i_1 k_1}{}_{j_1 k_2} \times T^{k_2k_3}{}_{j_2 j_3}
\end{displaymath}
for any $i_1, i_2, i_3, j_1, j_2, j_3$,
$1\le i_1(j_1)\le {\rm dim}\, V$,
$1\le i_2 (j_2)\le {\rm dim}\, V$, \linebreak $1\le i_3 (j_3)\le {\rm dim}\,
V$.
\vskip 6mm


{\bf Theorem.} \ {\it Any representation $\pi $ of algebra ${\cal B} $
$\pi : {\cal B} \to {\rm End}_k (V)$
 over field $k$, such that linear operator $t : V \oplus V \to V\oplus V$
with matrix $T=\left( \matrix { \pi (A) & \pi (B) \cr \pi (C) & \pi (D) }
\right)$
is nondegenerate, defines the representation of the  braid group $B_{\infty }$
over the field $k$.
%
  4
Any tensor $T\in V^{\otimes 2 } \otimes (V^* )^{\otimes 2}$, such
that relation {\rm (viii)} holds and matrix $(T^{uv}{}_{pq})$
is nondegenerate defines the representation of the braid group
$B_m$, for any $m$.   }
\vskip 4mm

The proof of the theorem is obvious if we set $V_i \simeq V$ for any $i\in N$.

In the case of theorem $d=1$ because $T_{i+1} =T_i$
(see remark 1.2.).
\vskip 4mm

{\bf Definition.} \ {\it We shall call an algebra ${\cal B}$
the braid algebra
%
  5
 ot type (2,1) of period 1 (or simply the braid algebra).
We shall call relation {\rm (viii)}  the braid equation of type (2,1) of period
1
(or simply the braid equation).  }
\vskip 4mm

Note that some decisions of the braid equation can be obtained via
quantum groups and Yang-Baxter equation, see [Dr], [T], [K].


\vskip 6mm

\noindent
2.3.  Let \ $V_1, V_2$ \ be two vector spaces over
the field $k$, $A_1, D_2 : V_1 \to V_1$,  $B_1, C_2 : V_2\to V_1 $,
$C_1, B_2 : V_1 \to V_2$,
$D_1, A_2 : V_2 \to V_2$ linear operators such that following
relations hold:
\noindent
$$
\matrix{
({\rm i}) &A_1^2 +B_1 A_2 C_1 =A_1;\hfill  &
({\rm vii}) &A_2 C_1 -C_1 A_1 =D_1 A_2 C_1;\hfill \cr
({\rm ii}) &A_2^2  + B_2 A_1 C_2 =A_2,\hfill &
({\rm viii}) &A_1 C_2 -C_2 A_2 =D_2 A_1 C_2,\hfill\cr
({\rm iii}) &D_2^2 +C_2 D_1B_2 =D_2;\hfill  &
({\rm ix}) &B_1A_2 -A_1 B_1 =B_1A_2 D_1;\hfill \cr
({\rm iv}) &D_1^2 +C_1 D_2B_1=D_1,\hfill    &
({\rm x}) &B_2 A_1 -A_2B_2 =B_2 A_1D_2 ,\hfill\cr
({\rm v}) &D_1B_2 -B_2D_2 =A_2 D_1 B_2;\hfill  &
({\rm xi}) &C_2 D_1 -D_2 C_2 =C_2 D_1 A_2;\hfill  \cr
({\rm vi}) &C_1D_2 -D_1C_1 =C_1 D_2 A_1,\hfill &
({\rm xii}) &D_2 B_1 -B_1 D_1 =A_1 D_2 B_1.\hfill \cr}
$$
$$
\matrix{
\hfill ({\rm viii}) & C_1B_1 -B_2 C_2 =A_2D_1A_2 -D_1A_2 D_1;\hfill  \cr
\hfill ({\rm xiv}) & C_2B_2 -B_1C_1 =A_1 D_2 A_1 - D_2 A_1 D_2,\hfill \cr}
$$
Let $T(1) \in (V_1 \otimes  V_2 )\otimes (V_1 \otimes V_2) ^*$,
 $T(2) \in (V_2 \otimes  V_1 )\otimes (V_2 \otimes V_1) ^*$ be two tensors,
such that for any $i_1, i_2, i_3, j_1, j_2, j_3$,
$1\le  i_1 (j_1) \le {\rm dim}\,V_1$,
$1\le  i_2 (j_2) \le {\rm dim}\,V_2$,
$1\le  i_3 (j_3) \le {\rm dim}\,V_1$, the following relations hold:

\begin{displaymath}
({\rm xv})
{\mathop {\sum} \limits_ {k_1, k_2, k_3}}
T(1) ^{i_1 i_2}{}_{k_1 k_3} \times T(2) ^{k_3 i_3}{}_{k_2 j_3}
\times T(1) ^{k_1 k_2}{}_{j_1 j_2} =
\end{displaymath}
$$
={\mathop {\sum} \limits_ {k_1, k_2, k_3}}
T(2) ^{i_2 i_3}{}_{k_1 k_3} \times T(1) ^{i_1 k_1}{}_{j_1 k_2}
\times T(2) ^{k_2 k_3}{}_{j_2 j_3}
$$
$$
(1\le k_1 \le {\rm dim} \,V_1,\quad 1\le k_2 \le {\rm dim} \,V_2,\quad
(1\le k_3 \le {\rm dim} \,V_1 ;)
$$

\begin{displaymath}
({\rm xvi})
{\mathop {\sum} \limits_ {k_1, k_2, k_3}}
T(2) ^{i_2 i_3}{}_{k_1 k_3} \times T(1) ^{k_3 i_1}{}_{k_2 j_1}
\times T(2) ^{k_1 k_2}{}_{j_2 j_3} =
\end{displaymath}
$$
={\mathop {\sum} \limits_ {k_1, k_2, k_3}}
T(1) ^{i_3 i_1}{}_{k_1 k_3} \times T(2) ^{i_2 k_1}{}_{j_2 k_2}
\times T(1) ^{k_2 k_3}{}_{j_3 j_1}
$$
\begin{displaymath}
(1\le k_1 \le {\rm dim} \,V_2,\quad 1\le k_2 \le {\rm dim} \,V_1,\quad
(1\le k_3 \le {\rm dim} \,V_2 ; )
\end{displaymath}
\vskip 6mm

{\bf Theorem.} \ {\it If linear operators $T_1 : V_1 \oplus V_2 \to V_1 \oplus
V_2 $ with
mat\-rix $T_1= \left ( \matrix{ A_1 & B_1 \cr C_1 & D_1 }\right) $,
$T_2 : V_2 \oplus V_1 \to V_2 \oplus V_1$ with  matrix
 $T_2= \left ( \matrix{ A_2 & B_2 \cr C_2 & D_2 }\right) $ are nondegenerate,
and relations {\rm (i)-(xiv) } hold then
linear operators $A_1, A_2, B_1, B_2, C_1, C_2, D_1, D_2$ define
representation of the braid  group $B_{\infty}$.

If matrices $(T(1)^{i_1 i_2}{}_{j_1j_2} )$,
 $(T(2)^{i_1 i_2}{}_{j_1j_2} )$ are nondegenerate and relations {\rm
(xv)-(xvi)}
hold then tensors $T(1)$ and $T(2)$ define representation
of the braid group $B_{m}$, for any $m\in N$.   }
\vskip 3mm

The proof of the theorem is obvious if we set $V_i =V_1$ for odd i and
$V_i =V_2 $ for even $i$.
\vskip 4mm

In the case of the theorem $d=2$, because $T_{i+2} =T_i$
(see remark 1.2.). It is obvious that analogous theorem can
be formulated for any $d$.
We should call the relations (i)-(xix) the braid algebra of type (2,1)
of period 2 and relations (xv)-(xvi) the braid equation of type (2,1)
of period 2.
\vskip 6mm

2.4. {\bf Definition.} \ {\it Let $k [x,x^{-1}]$ be the ring
of Laurent polynomials
over field $k$ and ${\cal B}^{\triangle} =k[x,x^{-1}] [y]$ is the
commutative algebra with the relation  $y^2 +xy =y$. We shall
call this algebra triangle braid algebra.  }
\vskip 6mm

{\bf Theorem.} \ {\it Suppose $\tilde t_i, i\in N$ (see 1.2 and 2.1 )
satisfies the conditions of the theorem  2.1   and $A_i =0$ $(D_i =0)$
for any
\linebreak $i\in N$. Then $V_i \simeq V_j \simeq V$ for any $i,j \in N$
and some vector space $V$,
\linebreak $W\simeq {\mathop {\oplus } \limits_{s\in N}}  V$
and the representation $\pi$, defined according to the theorem 2.1,
is equivalent to the representation $\pi '$ with $A_i' =0 (D_i' =0)$,
$C_i' =1_V$, $ B_i' =B$,
$D_i' =D (A_i' =A)$ for any $i$, where $D(A) :V \to V$,
$B: V\to V$ linear operators, and $B$ is nondegenerate (in particular, it is
the
case of theorem 2.2.). There exists one--to--one correspondence between
such representations $\pi '$ of  the braid group $B_{\infty}$
and the representations of the algebra ${\cal B}^{\triangle}$.    }
\vskip 8mm

\noindent

{\bf Proof. } \ Matrix of linear operator $\tilde t_i$ on the vector space $V_i
\oplus V_{i+1}$
is nondegenerate and it has the triangle form
$\left(\matrix{ 0 &B_i\cr C_i & D_i }\right)$
$\left(\left( \matrix{ A_i & B_i \cr C_i & 0 }\right)\right)$.
\vskip 2mm
\noindent
Hence
linear operators $C_i, B_i$ are nondegenerate and $V_i \simeq V_{i+1}$.
We can choose basis in every $V_i, i\in N$ such that matrix for any
$C_i , i\in N$ is $I$. After substituting $C_i =I, A_i =0 (D_i=0)$
into the relations 2.1. (i)--(vii) they come to the following form:
\begin{displaymath}
\matrix {BD =DB\hfill & (BA = AB),\hfill\cr
D^2 +BD =D\hfill & (A^2 +AB =A).\hfill\cr}
\end{displaymath}

Taking to account that $B$ must to be nondegenerate linear operator,
we obtain the statement of the theorem $(y\to D(A), x\to B)$.

\noindent

{\bf Remark.} \ {\it If the representation of the algebra ${\cal
B}^{\triangle}$ is direct
sum of some representations then corresponding representation of the
group $B_{\infty}$
 is the direct sum of the corresponding representations
of the  braid group $B_{\infty}$.   }
\vskip 4mm

The  equation $y^2 +xy =y$ is equivalent to the
equation $y(y+x -1) =0$. So it has at least two obvious solutions
$y=0$ and $y=1-x$. Thus for any nondegenerate linear operator
$B$ on the vector space $V$ over field $k$
we have two representations $f_1 : {\cal B}^{\triangle }\to {\rm End}_k V$
and $f_2 : {\cal B}^{\triangle }\to {\rm End}_k V$, where
$f_1 (x) =B$, $f_1(y) =0$, $f_2(x)=B$, $f_2(y) =I-B$.
So we have following cosequence of the theorem (we also took into account the
remark
above ) that gives three series I, II, III of the representations of the braid
group $B_{\infty}$:
\vskip 6mm

{\bf Consequence.} \ Let $B$ be linear nondegenerate indecomposable
operator $V\to V$ for some vector space $V$
 ( so if $k$ is algebraically closed field
and ${\rm dim} \, V< \infty $, then $B $ is Jordan box over field $k$
with nonzero eigenvalue).
\noindent
There are representations of the items 1.2.
and 2.1. with:

\noindent
\begin{displaymath}
\matrix{
({\rm i})\   ({\rm series I}) A_i =0, D_i =0, C_i =I, B_i =B \ \mbox{for any }
\ i\in N;\hfill\cr\cr
({\rm ii})\  ({\rm series II}) A_i =0,  C_i =I, B_i =B, D_i =I-B \ \mbox{for
any } \ i\in N;\hfill\cr\cr
({\rm iii})\ ({\rm series III}) A_i =1-B,  B_i =B, C_i =I, D_i =0 \ \mbox{for
any } \ i\in N;\hfill\cr \cr}
\end{displaymath}

\noindent

\noindent

\noindent
\vskip 8mm

\noindent
2.5. \ There are other more veiled representations of the commutative algebra
${\cal B}^{\triangle}$.
\vskip 4mm

{\bf Theorem. } \ {\it There is one-to-one correspondence between
finite-dimensional
indecomposable representations over field $k$ of the commutative algebra
${\cal B}^{\triangle}$, different from representations $f_1$ and $f_2$ (see
2.4.),
and the finite-dimensional
indecomposable representations over the field $k$ of the ordered pairs linear
operators
$(\pi_1, \pi_2)$ with relations
$\pi_1\pi_2 =0$ and $\pi_2 \pi_1 =0$ such that $\pi_1 \ne 0, \pi_2 \ne 0$
(corresponding algebra is, obviously, commutative algebra
$k[x,y]$ with the only relation   }
$xy =0)$.
\vskip 4mm

{\bf Proof.} \ Let  $f: {\cal B}^{\triangle} \to {\rm End}_k V $ be
indecomposable representation of the algebra $B^{\triangle}$,
$f\ne f_1, f\ne f_2$. Then $f(y) (f(y) +f(x) -I)=
(f(y) +f(x) -I) f(y)=0$. Let $\pi_1 =f(y), \pi_2 =f(y) +f(x) -I$.
Then $\pi_1 \pi_2 =\pi_2 \pi_1 =0$,
$\pi_1 \ne 0, \pi_2 \ne 0$
 and it easy to verify that
this representation of the pair of linear operators is indecomposable.
Conversely, let $(\pi_1 ,\pi_2)$ be the indecomposable representation of the
pair of
linear operators in vector space $V$, ${\rm dim} \, V<\infty $ such that
$\pi_1 \pi_2 =\pi_2 \pi_1 =0$, $\pi_1 \ne 0, \pi_2 \ne 0$.
If $\pi_1 (\pi_2 )$  is not nilpotent operator then $\pi_2 =0 (\pi_1 =0)$,
see [NRSB]. Hence $\pi_1$ and $\pi_2 $ are nilpotent operators.
Put $f(y) =\pi_1 , f(x) =\pi_2 -\pi_1 +I$.
As $\pi_2 $ and $\pi_1$ commutate one with another, linear operator
$\pi_2 -\pi_1 $ is nilpotent too and hence linear operator
$f(x)= \pi_2 -\pi_1 +I$ is invertible. Theorem is proved.
\vskip 4mm

{\bf Consequence.} \ {\it Every finite dimensional indecomposable
representation of the pair operators $\pi_1 ,\pi_2$ with relations $\pi_1 \pi_2
=
\pi_2 \pi_1 =0$ such that $\pi_1 \ne 0, \pi_2 \ne 0$ defines
representation of the braid group $B_{\infty}$.  }
\vskip 4mm

{\bf Proof.} \ It follows from the theorem and the theorem 2.4.
\vskip 4mm

\noindent
2.6. \ All indecomposable representations of the pair operators
$\pi_1, \pi_2 $ such that $\pi_1 \pi_2 =\pi_2 \pi_1 =0, \pi_1 \ne 0, \pi_2 \ne
0$
 have been found in [GP] for the case, when $k$ is algebraically closed field,
and in [NRSB] for the case of arbitrary field $k$.
\vskip 4mm

{\bf Theorem} [NRSB]. \ Every indecomposable finite-dimensional
representation of the pair operators $\pi_1, \pi_2$
with relations
$\pi_1 \pi_2 =\pi_2 \pi_1 =0, \pi_1 \ne 0, \pi_2 \ne 0$
has one of the two types:

({\rm i}) vector space $V$ of type $I$ is defined by sequence of pairs not
negative
integers $(k_1 ,l_1),...,(k_n,l_n)$, where $k_{\alpha} > 0$
when $\alpha > 1, l_{\beta} > 0$ when $\beta < n$.
Vector space $V$ is built like $1+ \sum (k_{\alpha} +l_{\alpha})$-dimensional
space
with basis of  vectors $v_{\alpha} \pi_1^k, v_{\alpha} \pi_2^l, v_n
\pi_2^{l_n}$,
where $k=0,...,k_{\alpha}, l=1,...,l_{\alpha}-1$,  $\alpha =1,...,n$.
Equality $\pi_1 \pi_2 =\pi_2 \pi_1 =0$,  form
of the basis vectors and relations
\begin{displaymath}
v_1 \pi_1^{k_1 +1} =0, \ v_{\alpha +1} \pi_1 ^{k_{\alpha +1}} =
v_{\alpha }\pi_2 ^{l_{\alpha}} \ (\alpha =1,...,n-1),
\ \ v_n \pi_2 ^{l_n +1} =0
\end{displaymath}
completely define the action of operators $\pi_1 $ and $\pi_2 $ on the
vectors of the basis;

(ii) vector space $V$ of type II is defined by

a) sequence $(r_1,s_1),...,(r_q,s_q )$ of pairs of natural numbers that is
nonperiodic . It means
that if this sequence coincide with sequence
$(r_1,s_1,),...,(r_t,s_t), (r_1,s_1),...,(r_t,s_t)$ then $t=q$.

Sequences that differ on cyclic substitution are equivalent,

b) polynomial $f(x) =x^t -b_1 x^{t-1} -\cdots - b_t (b_t \ne 0)$
over the  field $k$ that is the power of the indecomposable polynomial
over the field $k$.

Vector space $V$ of type II has basis of the vectors
$v_{\mu \nu}\pi_1^r, v_{\mu\nu}\pi_2^s$, where $r=0,...,r_{\mu}$,
$s=1,...,s_{\mu -1}$, $\mu=1,...,q$, $\nu =1,...,t$.
Operators $\pi_1$ and $\pi_2$ are completely defined by following relations:
\begin{displaymath}
\matrix{
v_{11}\pi_1^{r_1} =(v_{q1} b_1 +\cdots + v_{qt} b_t )\pi_2 ^{s_q},\hfill\cr
v_{1,\nu +1} \pi_1^{r_1} = v_{q\nu } \pi_2^{s_q} \ (\nu =1,...,t-1),\hfill\cr
v_{\mu +1, \nu } \pi_1^{r_{\mu +1}} = v_{\mu \nu }\pi_2 ^{s_{\mu}} \
(\mu =1,...,q-1, \nu =1,...,t-1).\hfill\cr  }
\end{displaymath}
\vskip 4mm

{\bf Consequence.} \ {\it There are series IV and V of the representations of
the
braid group $B_{\infty}$ originating from the representations of pair
operators $(\pi_1 , \pi_2 )$ of type I and II, respectively. }
\vskip 4mm

{\bf Proof.} \ It follows from the consequence 2.5. and the theorem.
\vskip 3mm

Theorem 2.5 and the theorem [NRSB] give the complete description of the
finite-dimensional
representations of the triangle braid algebra ${\cal B}^{\triangle}$.

In [NRSB] authors note the fact that the description
of Lorentz  group representations and the classification of finite
$p$-groups with abelian subgroup of index $p$ are connected with the
same problem of description mutually annihilating operators is remarkable.
Now we see that we have to mention  that this problem is connected with
braid group representations too.
\vskip 4mm

\noindent
2.7. {\bf Theorem.} \ {\it Let $\pi :{\cal B} \to {\rm End}_k V$ be the
representation of the braid
algebra ${\cal B} $ over field $k$, such that
$\pi (C)=I$. Then the corresponding
representation of the braid group $B_{\infty }$ (see theorem 2.2.),
in such case, when operator on the vector space $V\oplus V$ with the
matrix $\left( \matrix{ \pi (A) & \pi (B)\cr I       & \pi (D)}\right )$
is nondegenerate, may be obtained of some representation $\pi '$
of algebra
\ $\widetilde {{\cal B}}$ \,over field $k$ with generators
$A, B, D$ and the following defining relations:
\begin{displaymath}
\matrix{
({\rm i}) & A^2 +BA =A,\hfill\cr
({\rm ii}) & D^2 +DB =D,\hfill\cr
({\rm iii}) & BA -AB =BAD,\hfill\cr
({\rm iv}) & DB -BD =ADB,\hfill\cr
({\rm  v}) & DA=0.\hfill\cr }
\end{displaymath}

Conversely, any representation $\pi$ of algebra $\widetilde {{\cal B}}$
$\pi :{\cal B} \to {\rm End}_k V$,
 such that
operator on vector space $V$ with matrix
$\left( \matrix { \pi (A) & \pi (B) \cr
I & \pi (D) }\right) $ is nondegenerate, gives the representation
of the braid group $B_{\infty}$.     }
\vskip 4mm

{\bf Proof.} \ It is obvious.
\vskip 4mm

{\bf Definition.} \ We shall call the algebra
$\widetilde {{\cal B}}$ simplified braid algebra.

\noindent
\vskip 6mm

\noindent
2.8. {\bf Theorem.} \ {\it Quotient algebra of the braid algebra
${\cal B}$
on it commutant  is isomorphic to the commutative algebra
$ {\cal B}^{{\rm com}} =k[x,y,z,t]$ with following difining  relations:
\begin{displaymath}
\matrix{ ({\rm i}) & x^2 +xyz =x,\hfill\cr
         ({\rm ii}) & t^2 +tyz =t,\hfill\cr
        ({\rm iii}) & xzt =0,\hfill\cr
({\rm iv}) & xyt =0,\hfill\cr}
\end{displaymath}

Any representation $\pi $ of the algebra $ {\cal B}^{{\rm com}} $
$\pi :{\cal B}^{{\rm com}} \to {\rm End}_k V$,
 over the field $k$ such that operator on the vector space
$V\oplus V$ with matrix $\left( \matrix {\pi (x) & \pi (y) \cr
\pi (z) & \pi (t)}\right)$ is nondegenerate defines the representation of the
braid group $B_{\infty }$.    }
\vskip 4mm

{\bf Proof.} \ It is obvious (here $A\to x, B\to y, C\to z, D\to t)$.
\vskip 4mm

{\bf Definition.} \ We shall call algebra $ {\cal B}^{{\rm com}}$ commutative
braid
algebra.
\vskip 4mm

{\bf Remark.} \ We note that if the representation $\pi $ of the braid group
$B_{\infty }$
is obtained from some representation of the commutative braid algebra, it
doesn\'{}t
mean, of course, that $\pi (B_{\infty })$ is commutative.
\vskip 4mm

\noindent
2.9. {\bf Theorem.} \ {\it Quotient algebra of the simplified braid algebra
 $\widetilde {{\cal B}}$  on its commutant is
isomorphic to the commutative algebra
$\widetilde {{\cal B}}^{{\rm com}} =k[x,y,t]$ with the following defining
relations:
\begin{displaymath}
\matrix{ ({\rm i}) & x^2 +xy =x,\hfill\cr
({\rm ii}) & t^2 +ty =t,\hfill\cr
({\rm iii}) & xt =0.\hfill\cr}
\end{displaymath}

Any representation $\pi $ of the algebra $\widetilde {{\cal B}}^{{\rm com}} $
$\pi :\widetilde{{\cal B}}^{{\rm com}} \to {\rm End}_k V$,

in the vector space $V$ over the field $k$, such that operator on the
vector space $V\oplus V$ with matrix $\left (\matrix {
\pi (x) & \pi (y)\cr I & \pi (t) }\right)$ is nondegenerate, defines the
representation of the braid group $B_{\infty }$.    }
\vskip 4mm

{\bf Proof.} \ It is obvious (here $A\to x, B\to y, D\to t)$.
\vskip 4mm

{\bf Definition.} \ We shall call algebra $\widetilde {{\cal B}}^{{\rm com}} $
simplified commutative braid algebra.
\vskip 6mm

\noindent

{\bf Remark.} \ As items 2.5 and 2.7 show any representation of  the braid
group
$B_{\infty}$, obtained of some finite-dimensional representation of algebra
$ {{\cal B}}^{{\rm com}} $, may be obtained from  some finite-dimensional
representation of algebra $\widetilde {{\cal B}}^{{\rm com}} $ .
\vskip 4mm

2.10. {\bf Remark.} \ If the representation of the algebra
$\widetilde {{\cal B}}^{{\rm com}} $  is direct sum of some representations,
then corresponding representation of the braid group $B_{\infty }$ is
the direct sum of the corresponding representations of the braid group
$B_{\infty }$.
\vskip 3mm

Let $R_1$ be the  quotient algebra of algebra $ \widetilde {{\cal B}}^{{\rm
com}} $
on the  relation $t=\alpha x, \alpha \in k, \alpha \ne 0$. Then $R_1$ is
commutative algebra $k[x,y]$ with the following relations:
\begin{displaymath}
\matrix { ({\rm i}) & x^2 =0,\hfill\cr
({\rm ii}) & x(y-1) =0.\hfill\cr}
\end{displaymath}

The relation (ii) holds automatically if $y=1 +\beta x$, for  any $\beta \in
k$.
So we have the commutative algebra $R_2 =k[x]$ with the only relation
$x^2 =0$. Its indecomposable representation is defined by the Jordan
$2\times 2$ box with zero eigenvalue $x\to \left( \matrix{ 0 & 1\cr
0 & 0 }\right)$. So we have representation $\pi_{\alpha, \beta }$
of algebra $ \widetilde {{\cal B}}^{{\rm com}} $ in the vector space $V$
with basis $e_1$, $e_2$ such that
$\pi_{\alpha ,\beta } (x)= \left( \matrix{  0 & 1\cr 0 & 0}\right) $,
$\pi_{\alpha ,\beta } (y)= \left( \matrix{  1 & \beta \cr 0 & 1}\right) $,
$\pi_{\alpha ,\beta } (t)= \left( \matrix{  0 & \alpha \cr 0 & 0}\right) $.
It is easy to verify that $\pi_{\alpha ,\beta }$ is indecomposable.
Operator on the vector space $V\oplus V$ with matrix
$\left( \matrix{  \pi_{\alpha ,\beta } (x) & \pi_{\alpha ,\beta } (y)\cr
I & \pi_{\alpha ,\beta } (t) } \right) =
\left( \matrix{  0&1 &1 &\beta \cr
0&0&0&1\cr 1&0&0& \alpha\cr  0&1&0&0\cr }\right) $
is nondegenerate, hence, according to theorem 2.9,
$\pi_{\alpha ,\beta }$ defines the representation of the braid group
$B_{\infty }$.
Matrices of the operators $\tilde t_i, \tilde t_{i+1}$ in the vector space
$W ={\mathop {\oplus } \limits_{s\in N}} V_s (V_s \simeq V$ for any $s\in N)$
in the basis
$e_1^{(s)}, e_2^{(s)}, s=1,2,...$ (not indicated matrix elements are 1 if they
are
diagonal, and $0$ if not) have the next form:
\begin{displaymath}
\tilde t_i =\left( \matrix{
e_1^{(i)} & e_2^{(i)} & e_1^{(i+1)} & e_2^{(i+1)} & e_1^{(i+2)} & e_2^{(i+2)}
\cr
0 & 1 & 1 & \beta & 0 & 0 \cr
0 & 0&0&1&0&0\cr
1&0&0&\alpha &0&0\cr
0&1&0&0&0&0\cr
0&0&0&0&1&0\cr
0&0&0&0&0&1  }\right) ;
\end{displaymath}
\begin{displaymath}
\tilde t_{i+1} = \left( \matrix {
e_1^{(i)} & e_2^{(i)} & e_1^{(i+1)} & e_2^{(i+1)} & e_1^{(i+2)} & e_2^{(i+2)}
\cr
1& 0& 0& 0& 0& 0\cr
0& 1& 0& 0& 0& 0\cr
0& 0& 0& 1& 1& \beta \cr
0& 0& 0& 0& 0& 1\cr
0& 0& 1& 0& 0 & \alpha \cr
0& 0& 0& 1& 0& 0 } \right).
\end{displaymath}

One  can easily verify that
$\tilde t_i \tilde t_{i+1} \tilde t_i =
\tilde t_{i+1} \tilde t_i \tilde t_{i+1}$.
This series of the braid group representations we shall call series VI.

The next theorem generalize the above example of the representation
of algebra $\widetilde{{\cal B}}^{{\rm com}}$.
\vskip 6mm

{\bf Theorem.} \ {\it  Let $R$ be a finite-dimensional algebra over the
field $k$ and $J$ is an ideal in $R$, such that $J^2 =0$.
 For the elements $a,b, c \in J$ put $\varphi (x)=a$, $\varphi (t)=b$,
$\varphi (y)=1+c$. Then $ \varphi $  determines homomorphism of algebras
over the field $k$  $ \varphi  : \tilde {\cal B}^{{\rm com}} \to R$.
Put $\pi =\psi  \varphi $, where $\psi$ is a regular representation
of algebra $R$ over the field $k$. Then $\pi$ defines  finite-dimensional
linear representation of algebra  $\tilde {\cal B}^{{\rm com}}$,
such that matrix $\left(\matrix{
\pi (x) &  \pi (y)\cr I  &  \pi (t) }\right)$ is nondegenerate
and, hence, the representation of the group $B_{\infty}$.  }
\vskip 5mm

{\bf Proof.} It is obvious.
\vskip 4mm

{\bf Definition.} \ {\it The algebra over the field $k$ (the group) is called
wild if the problem of classification of its linear  representations contains
the problem of pair matrices (that is the problem of classification for any
pair matrices $A, B \in {\rm Mat}_n (k)$, $ n> 1$ the orbits
$XAX^{-1}$, $XBX^{-1}$, $X\in GL_n (k)$), see also [D1].  }
\vskip 4mm

Note that if the algebra is wild it has, roughly speaking, "a huge amount"
of indecomposible representations.
\vskip 6mm

{\bf Consequence.} \ {\it Algebra $\tilde {\cal B}^{{\rm com}}$ (and,
hence, ${\cal B}$) is wild.   }
\vskip 4mm

{\bf Proof.} \ Let
\noindent
$$
a=\left(\matrix{0  & A \cr 0 & 0}\right)_{2n\times 2n},
b=\left(\matrix{ 0  & B \cr 0  &  0}\right)_{2n\times 2n},
c=\left(\matrix { 0 & C \cr 0  &  0}\right)_{2n\times 2n},
$$
where $A,B,C \in {\rm Mat}_n (k)$, $R=k< a,b,c > \subset {\rm Mat}_{2n} (k)$,
$J=RaR +RbR +RcR$. Then we appear in the conditions of the theorem.
But the problem of classification
of the triple of matrices $a,b,c+1 \in {\rm Mat}_{2n} (k)$ is wild,
see [D2].
\vskip 4mm

{\bf Remark.} \ {\it The above consequence, it seems, is a very favourable for
the theory of links invariants.}
\vskip 4mm

Series I-VI, theorem and the consequence give, in certain sense, the decision
of the problem 25 of Joan Birman (Problem 25. Finite
representations of the braid groups), see [M].

If $k$ is a finite field we get from the series I-VI and the theorem the
finite quotiens of $      B _n$ (in fact, most of obtained
representations can be defined, when $k$ is a ring, in particular,
when $k$ is  a  finite ring,for instance $k={     Z} /m {     Z},\,{m\in N}$).
It gives a lot of decisions of the problem 25.1 (Problem 25.1.
Identify explicit interesting finite quotients of ${ B}_n$),
see [M].

Of course, the fact, that we have found simple quotient algebra of the
braid algebra ${\cal B}$ (namely the commutative braid algebra) mainly
contributed to the obtaining of the results of this item.
\noindent
2.11. \ There are other interesting quotient algebras of the braid algebra. For
instance, if we factorize the braid algebra on the  relation $AD=0$, we obtain
noncommutative algebra over the field $k$ with generators
$A, B, C, D$ and with the following relations:
\begin{displaymath}
\matrix{ ({\rm i})  & A^2 +BAC=A,\hfill\cr
({\rm ii}) & D^2 +CDB=D,\hfill\cr
({\rm iii}) & CB=BC,\hfill\cr
({\rm iv}) & AB=BA,\hfill\cr
({\rm v}) & AC-CA=DAC\hfill\cr
({\rm vi}) & CD-DC=CDA\hfill\cr
({\rm vii}) & AD=0.\hfill\cr  }
\end{displaymath}


\noindent
\vskip 8mm

\noindent
2.12. \ We may suppose that there are "mixed" series, "glued together"
of obtained series. We mean that possibly there is a representation of the
braid group $B_{\infty }$ such that
for some $n_1 \in N$ $\tilde t_i$ could be of the series $k_1$ for $i\in
[1,n_1]$,
then for  some $n_2 > n_1,t_i$, are of the "glue"
if $i\in [n_1 +1, n_2]$, then for some
$n_3 > n_2 $ $\tilde t_i$ could be of the series $k_2$ if $i\in [n_2 +1, n_3]$
and so on. It is interesting to investigate if this could happen and if
"yes", to describe "all kinds of glue"  that could "stick" series $k_1$ and
$k_2$.

\vskip 8mm

\noindent
2.13. \ The theory of the representation of the infinite group $G$ includes
 the description all nonizomorphic short exact sequences
\begin{displaymath}
0 \longrightarrow M\longrightarrow N\longrightarrow  P \longrightarrow 0
\end{displaymath}
of $G$ modules. In this paper we give only some examples of $B_{\infty}$
modules,
leaving the investigation of their submodules and quotient modules for
 further publications as well as  consideration on the base of the
same approach of Briescorn group\'{}s representations, bridge  algebra
representations and so on (we mean the theory of  the representations of
algebras and groups with transferable invariant defining relations, see remark
1.2.).

We leave also for further publication functional aspects of the approach of
this article,
that is, when $k=R$ or $k=C$ and $V_i$ in $W={\mathop {\oplus } \limits_{i\in
N}} V_i$ have
the structure of Hilbert or unitar spaces  and we consider unitary,
self-conjugate operators and so on.
\vskip 4mm

\noindent
2.14. \ It is easy to formulate the statement analoguos theorem 2.2 (that is,
for $d=1$)
for the case $(n,\lambda )=(3,1);$ $(n,\lambda )=(3,2 );$ $(n,\lambda )
=(4,1);$
$(n,\lambda )=(4,2)$ etc.
We shall call the corresponding tensor
the braid tensor of type $(n,\lambda)$ and
the corresponding algebra the
braid algebra of type $(n,\lambda )$ and denote it ${\cal B}_{n,\lambda }$.
Thus ${\cal B} ={\cal B}_{2,1}$.
Despite the big number of generators and defining relations,
some of the factorizations of these algebras are "not complex".
We note that there is an obvious map from the set of
representations of algebra ${\cal B}_{nm, \lambda m}$ for any $m\in N$ into the
set of the representations of algebra ${\cal B}_{n,\lambda }$.
Some of the representations of these algebras ${\cal B}_{n,\lambda }$
(now hypothetical ) generate representations of the braid group
$B_{\infty }$
that "detail" those obtained from the representations of the braid algebra
${\cal B}_{n_1,\lambda_1 }$
(not only in the case when $n=n_1m, \lambda =\lambda_1 m$ for some $m\in N$ )
in particular of the braid algebra ${\cal B} $. The amount of the possibilites
and combinatorics that appear here are impressive.


\noindent
\vskip 6mm

\noindent
2.15. \ "{\large {\sl Irresponsible dynamical system interpretation}}".

Suppose that all
representations of the braid group $B_{\infty}$ are the points of the
variety (or may be the projective or injective limit of varieties).
Map $F$ of this variety (into itself) can be described in the
following way: to every representation $\pi $ (the point of variety)
of the braid group $B_{\infty}$ we set point $F(\pi )$, namely
representation $F(\pi )$ such that
$F(\pi )(g_i)=\pi (g_{i +1})$ for any $i\in N$. Thus we have "the dynamical
 system "on the "variety" of all representations of the group $B_{\infty}$.
The obvious analog of the theorem 2.2 for the general case $(n,\lambda)$
and  $d> 1$ gives the
description of some
periodic orbits of such dynamical system, and the obvious analog of
theorem 2.2 for the algebras ${\rm B}_{n,\lambda}$ in case $d=1$
gives the description of some
fixed points of this dynamical system. If the sequence $\{F^n(\pi)\}_{n\ge 1}$
"converges" to some point $\pi '$, then $\pi '$ is a fixed point.
The interesting question is, can any fixed point (or more widely,
can any periodic orbit) of this dynamical system be obtained by the
approach utilized in  this paper?
\vskip 5mm

\noindent
2.16. \ For any $n\in N$ let $B_n\subset B_{\infty}$ be the group generated by
the elements $t_1,...,t_{n-1}$. It is called the braid group with $n$
strings [B]. Let $\pi$ be a finite-dimensional representation of the group
$B_n$. Algebra $k\! < \pi (t_1),..., \pi (t_{n-1}) >$ generated over
the field $k$ by operators $\pi (t_1),...,\pi (t_{n-1})$ has
the finite dimension over the field $k$. For the group algebra $kB_n$
of $B_n$ over the field $k$ let ideal $J_{\pi}$ be the kernel of the map
$kB_n \to k \!<\pi(t_1),...,\pi (t_{n-1})>$ $(t_i\to \pi (t_i))$.
Let $J$ be any ideal of the algebra $kB_n$, such that
$J\subset J_{\pi }$ and the quotient algebra $kB_n/J=H_{\pi,J,n}$ has
the finite dimension over the field $k$. Such algebras $H_{\pi , J,n}$
we call "Hecke alike" algebras generated by the representation $\pi $.
Among the linear representations of such algebra there is one that
defines representation $\pi $ of the braid group $B_n$,
namely $kB_n\to kB_n /J\to kB_n / J_{\pi }$. Thus,
algebras $H_{\pi ,J, n}$ "multiply" the representation $\pi $.
The most interesting is the case when $J$ is minimal ideal in $kB_n$
among those with the finite dimensional quotient algebra.
If we have finite-dimensional
linear representation $\pi $ of the group $B_n$ it is easy, as the
rule, to find out the relationships in the ideal $J_{\pi}$ that "make" the
dimension of the quotient algebra finite.
 Choosing these relationships and factorizing
on them, we obtain finite-dimensional "Hecke alike" algebra,
generated by $\pi $. If the representation $\pi $ "depends" of some parameter
$\lambda $ the algebra $H_{\pi, J,n}$ can often be "lifted" to the algebra
over the ring $k[\lambda , f(\lambda )^{-1}]$ $(f(\lambda ))$
polynimial over the field $k$), where $\lambda $ is undeterminate
(it sometimes becomes difficult to find that this algebra
originates from the representation $\pi$).
Algebras obtained in this way, it seems, will play
similar role
in the construction of invariants for knots and links, as it plays algebra
Hecke
$H(n,q)$ [J] (the most important problem is to find the "trace" with
special properties on these algebras).

Let us illustrate the  above  case when $\pi$ is from the series II,
$k= {\cal C}, B=q$, where $q\in {\cal C}, q\ne 0$ (it is the Burau
representation, see [B]).

We have $\tilde t_i =\left( \matrix{ 0&q\cr 1& 1-q }\right)$
(see 1.2, 2.1, 2.4(ii)). It is easy to verify  that
$({\rm i}) \ \ \tilde t^2_i -(1-q) \tilde t_i -q \cdot I=0$.
Representation $\pi $ defines the map: ${\cal C} B_n \to {\cal C}\! <\pi
(t_1,...,\pi(t_{n-1})>$. Let $J_{\pi }$ be the kernel of this map.
Relationship (i) belongs to $J_{\pi }$ (and alone "makes"
the dimension of the quotient algebra ${\cal C}B_n / J_{\pi}$ finite).
Let $J$ be the ideal in ${\cal C} B_n$ generated by relationship (i).
Then ${\cal C}B_n /J =H(q,n)$, where $H(q,n)$ is the  algebra
Hecke, ${\rm dim}_{{\cal C}} H(q,n)=n! <\infty$,
$({\rm dim}_{{\cal C}} {\cal C}\! < \pi (t_1),...,\pi(t_{n-1})> < n^2$,
hence $J \subset J_{\pi}$ and $J\ne  J_{\pi}$).
The representations of the algebra $H(q,n)$ (that is generated by $\pi $)
in particular "multiply" the representation $\pi $.
If $q$ is undeterminate, we can "lift" the algebra $H(q,n)$ to the
algebra over the ring ${\cal C}[q,q^{-1}]$ of the Laurent
polynomials. Algebra that appears plays an important role in the construction
of polynomial invariants of knots and links [J], [F--O].

We see, that algebra Hecke $H(q,n)$ is one of the  "Hecke alike"
algebras, that is generated by the representations of the braid group
$B_{\infty}$, that appears here.
We hope to devote our following publication to the

\noindent
 study and investigation  of these "Hecke alike" algebras (in
particular generated by representations from the series I--VI)
and invariants of links that possibly could happen here.
\vskip 5mm

\noindent
2.17. \  We indicate one case when the invariant of links
appears directly (and the "trace" is "isual" trace).

Let $\pi $ be the finite dimensional representation of the braid
algebra ${ \cal B}$ into ${\rm End}_k(V)$ ($V$ is the vector space over the
field $k=R$ or $k={\cal C})$, such that linear operator
$V\oplus V\to V\oplus  V$ with matrix $\left( \matrix{ \pi(A) & \pi(B) \cr
                                         \pi(C)  &  \pi (D)}\right) $ is
nondegenerate. Then there exist elements $A_1, B_1, C_1, D_1$ of the
algebra ${\cal B}$ such that the following equality holds
\begin{displaymath}
\left( \matrix{ \pi (A) & \pi (B) \cr
\pi (C) & \pi (D) }\right)
\left( \matrix{  \pi (A_1) & \pi (B_1) \cr
\pi (C_1) & \pi (D_1) }\right) =\left( \matrix{ I & 0\cr 0 & I}\right).
\end{displaymath}
Note that linear operators $\pi (A_1), \pi(B_1), \pi(C_1), \pi(D_1)$
on the vector space $V$ are uniquely determined by the representation $\pi$.
\vskip  4mm

{\bf Definition.} \ {\it We shall call the above representation $\pi$ of the
braid algebra
"simple" representation if the following conditions hold:

\begin{displaymath}
\matrix {
({\rm i}) \ {\rm tr }\,\pi (AX) ={\rm tr}\, \pi (X) + t q_1(X), \quad q_1(X)\in
{ Z} ;\hfill\cr\cr
({\rm ii}) \ {\rm tr }\,\pi (A_1X) ={\rm tr}\, \pi (X) + t q_2(X); \quad
q_2(X)\in { Z} \hfill }
\end{displaymath}
for some $t\in k, t\ne 0$ and for any $X\in {\cal B} $ one of the following
forms:
\begin{displaymath}
\matrix {
({\rm iii}) \ X  \ \mbox {is the product of some  elements }\hfill\cr
                \hfill  \mbox {of algebra}\ {\cal   B}, \mbox{ each equals} \ D
\mbox{ or} \ D_1,
                  \mbox{ or} \ X=1;\cr\cr
({\rm iv}) \ X=G_1S_1MS_2G_2, \ \mbox{where} \ G_1(G_2) \hfill }
\end{displaymath}
is 1 or the product of some elements of algebra $\cal B$ each equals $D$ or
$D_1$, $S_1 (S_2)$
equals $C$ or $C_1 (B$ or $B_1)$, $M$ is the product of some elements of
algebra
$\cal B$, each equals to one the elements $A, A_1$, $B, B_1$, $C, C_1$,
$D, D_1$, such that $n_1 +n_2 =n_3 +n_4$, where
$n_1, n_2, n_3, n_4$ are the numbers of elements $B, B_1$, $C, C_1$
respectively, in this product. }
\vskip 5mm

{\bf Theorem.} \ {\it Let $\pi$ be "simple" representation of the braid
algebra $\cal B$. For any oriented link $\hat \beta$, obtained of the braid
$\beta \in B_n$ (see [B]),
let $ I(\hat\beta) = e^{{\pi i \over t}
[ 2{\rm tr}\, \pi' (\beta) +(\exp \beta ) {\rm tr}\, [\pi (D_1) -\pi (D)]
- n {\rm tr}\,[\pi (D_1) +\pi (D) ]] }$
where $\pi'$ is the representation of the group $B_n$ into ${\rm End}_k
\left( {\mathop {\oplus } \limits_{i=1}^{n}} V \right)$
obtained of the representation $\pi$ (theorem 2.2),
$\exp \beta = \alpha_1+\cdots +\alpha_k$, if
$\beta =t_{i_1}^{\alpha_1} ...t_{i_k}^{\alpha_k}$. Then $I(\hat \beta )$
is the invariant of the oriented link $\hat \beta $.  }
\vskip 4mm

\rm
{\bf Scetch of the proof.} \ It is obvious that $I(\hat\beta )$ does not depend
of Markov move of type 1 (see [B]). Besides this, it does not
depend of Markov move of type 2 (see [B]) if we take into account,
that ${\rm tr}\,  \pi' (t_n \beta)= {\rm tr}\, (\beta) +
{\rm tr}\,  D +tq_1 (X_1)$, ${\rm tr}\,  \pi' (t^{-1}_n \beta )=
{\rm tr}\,  \beta +{\rm tr}\, D_1 +tq_2(X_2)$,
where $X_1(X_2)$ is the sum of some elements of algebra $\cal B$ each of the
form
that is indicated in definition, $q_1 (X_1) (q_2 (X_2))$ is the
corresponding sum of the values of the function $q_1(q_2)$ on these
summands.
Thus, according to Markov theorem (see [B]), $I(\hat\beta )$ is
invariant of the link $\hat\beta $.
\vskip 4mm

{\bf Remark.} \ {\it The theorem remains true if we shall demand the holding
of the relations {\rm (i), (ii) } for the $X$ of the form
{\rm (iii), (iv)}, and such that
$\psi (i) =n_1 (i) +n_2(i) -n_3(i) -n_4(i) \ge 0$ for any
$1\le i < l(X)$, where $n_1(i)$, $n_2(i)$, $n_3(i)$, $n_4(i)$
are the numbers of appeareces of generators $B, B_1, C,C_1$,
respectively, in monom $X$, from the first position of the monom to the
$i^{{\rm th}}$ position, $l(X)$ --- the length of the monom $X$.
See also 3.9, 3.10, 3.11.        }

\vskip 6mm

\noindent
{\Large\bf 3} \ {\large {\bf Some invariants  of oriented knots and links}}
\vskip 3mm

In this part of the paper  we indicate new methods of constructing invariants
of knots
and links (different from "classical" method, see  item 2.16 and [J]).
We find some invariants of oriented knots and links.
\vskip 4mm

\noindent
\rm 3.1. {\bf Remark.} \ {\it There exists well known geometric realization of
braid
group $B_n$, (see [B]). For geometric braids $\beta_1, \beta_2 \in B_n$,
$\beta_1 \cdot \beta_2$ is the result of concatenation of braids $\beta_1 ,
\beta_2$
$(\beta_2 \ $-- first). }

\vskip 3mm
\noindent
3.2. {\bf Definition.} \ {\it Let $G$ be any group. By the
$G$-braid of $n$ strings we shall understand any permutation $\pi $ on
the set $\{ 1,2,...,n\}$, with the map, that to any pair $(i,\pi (i))$
$i\in N, i\le n$ corresponds some element $g_i \in G$. }
\vskip 3mm

\rm
It is obvious what way we have to multiply two $G$-braids of $n$-strings and
that we obtain the group. Let denote it $B_n(G)$. It is obvious, that
$B_n(G) \subset B_{n+1} (G)$, so we can consider the group
$B_{\infty} (G) = {\mathop {\cup }\limits_ {n}} B_n(G)$.

\vskip 3mm
\noindent
3.3.  Let $a_i,b_i$, $i=1,\infty $  be the sequence of elements of the group
$G$, such that the following relations hold:
\begin{displaymath}
\matrix {({\rm i}) & a_i b_i  =b_{i+1} a_{i+1}  \hfill\cr } \qquad \mbox{for\
any }\
i\in N
\end{displaymath}
(these relations appeared first in [M]).

Let us correspond to any $t_i \in B_{\infty}$ the $G$-braid in
$B_{\infty}(G)$, such that $\pi(i)=i+1$,
$\pi (i+1)=i$, $\pi (k)=k$, if $k \not \in \{i, i+1\}$,
$g_i=a_i$, $g_{i+1} =b_i$, $g_k=1$, if $k \not\in \{i, i+1\}$.
We shall denote this $G$-braid $\tilde t_i$.
\vskip 6mm

\rm
{\bf Theorem.} \ {\it The map $t_i \to \tilde t_i$ defines
homomorphism $B_{\infty} \to B_{\infty} (G)$.}
\vskip 3mm

\rm
{\bf Proof. } \ It is easy to verify that it is equivalent to the
relations (i).
\vskip 3mm

The image of the braid $\beta \in B_{\infty}$ under this
homomorfism we shall denote by $\tilde \beta (\tilde \beta \in B_{\infty}
(G))$.
We see, that to every string of the braid $\beta$ corresponds some element
of the group $G$.
\vskip 3mm

\noindent
3.4. Given $\beta \in B_n$, one can correspond a braid in the space and may
form
an oriented link $\hat \beta $ by identifying the point $j$ at the top and
the point $j$
in the bottom ($j=1,...,n$) (see [B]). Any oriented link
arises in this way, and the question of when two
braids give rise to isotopic links is answered by
\vskip 6mm

{\bf Theorem.} (A.A. Markov, see [B]). \ {\it Let $B$ be the disjoint union
\begin{displaymath}
B={\mathop {\bigsqcup }\limits_ {n\ge 1}} B_n
\quad (B_1 =\{ 1\} ).
\end{displaymath}
Define the equivalence relation on $B$ as the one which is generated by the
relations
\begin{displaymath}
\matrix { ({\rm i}) & \mbox{If} \ \beta \ \mbox {and} \ \gamma \in B_n,
              \ \mbox{then} \ \beta \equiv \gamma \beta \gamma^{-1}, \hfill\cr
({\rm ii}) & \mbox {If} \ \beta \in B_n, \ \mbox {then} \ \beta \equiv t_n\beta
                                    \ \, \mbox{and }           \hfill\cr
& \beta \equiv t_n^{-1}\beta \ (\mbox{where} \ t_n^{\pm 1} \beta \in
B_{n+1}).\hfill\cr}
\end{displaymath}
Then two braids $\beta, \beta '$ give rise to isotopic links if and only
if $\beta \equiv \beta '$.
\vskip 3mm

\rm
\noindent
3.5. Next simple theorem gives the key to this part of the paper.
\vskip 6mm

{\bf Theorem.}  \ { \it Let $G$ be any group and $a_i,b_i$ $i\in N$
be any sequence of its elements that satisfy the relations 1.3 (i). For $\beta
\in B_n$
let us go from any point $j\in N$ (in the top or in the bottom) of some knot
in $\hat \beta $ in the direction of this knot  multiplying subseguently
corresponding strings in $\tilde \beta \in B_n(G)$,
that are the elements of $G$, (the stick of the point $j$ in the top with
the point $j$ in the bottom corresponds to $1\in G$) until we return to the
starting point. Let us correspond to the every knot in $\hat \beta $ the
conjugacy
class of the result of this multiplication. Then if  $b_i =a_i^{-1}$,
the family of this cojugacy classes of the group $G$
(their number is the number of knots in $\hat \beta $)
 is  invariant of the link $\hat \beta$. }
\vskip 3mm

\rm
{\bf Proof. } \ Let $\beta_2 \in B_n$ be obtained from $\beta_1 \in B_n$ by
Markov
move 3.4 (i). Then $\beta_2 =\gamma \beta_1 \gamma^{-1}$ for some
$\gamma \in B_n$ and, hence, $\tilde \beta_2 =\tilde \gamma \tilde \beta_1
\tilde \gamma^{-1}$.

Let us go from any point of some knot in $\hat \beta_2$ in the direction of
orientation
of link $\beta_2$ multiplying subsequently elements of group $G$
corresponding to the strings, that we are passing through, until we return to
the starting point.
It is not difficult to see, taking into account the equality
$\tilde \beta_2 =\tilde \gamma \tilde \beta_1 \tilde \gamma^{-1}$, that
the conjugacy class of the result of multiplication is equal to the
conjugacy class for the corresponding knot in
$\hat\beta_1$. Thus the families of conjugasy classes so obtained
for link $\hat\beta_1$ and for link
$\hat\beta_2$, coincide.

Let $\beta_1 \in B_n$ and $\beta_2 =t_n^{\pm} \beta_1$, i.e. $\beta_2$ is
obtained of $\beta_1$ by Markov move 3.4 (ii)).
Then, evidently, the only one conjugacy class for $\hat\beta_2$ could changes
(according to those in $\hat\beta_1$), namely the conjugacy class that
corresponds
to the  knot passing through the $(n+1)^{{\rm th}}$ point in the bottom. But
 it does not  change because the only new multiplier that
appears in corresponding product is element $b_n a_n=1$.
Thus the family of conjugacy classes does not depend on
Markov moves, that proves
the theorem.
\vskip 6mm

{\bf Consequence.} \ \it For the field $k$ and any $m\in N$ let
$K=k (\{ a_{ij}^{(s)} \}$) be the field of
the rational functions of
the undeterminates $\{ a_{ij}^{(s)}\}, i,j =1,...,m$;
$s\in N$ over the field $k$ (every function depends of
the finite number of undeterminates). For the sequence
$a_s=( a_{ij}^{(s)}) $, $b_s=(a_{ij}^{(s)})^{-1}$, $s=1,\infty$ of the
elements of $GL_m (K)$ and
every link $\hat \beta  (\beta \in B_n$ for some $n\in N$) consider
the family of characteristic polynomials of any integer power $t$ of
representatives of
the conjugacy classes, indicated in the theorem. Then this family of
polynomials of $X$ in $k[ \{a_{ij}^{(s)} \}$,
${1\over {\rm det}\,(a_{ij}^{(s)}) },X] $ is invariant of the link.
\vskip 4mm

\rm
{\bf Proof. } \ It is obvious.
\vskip 4mm

{\bf Remark.} \ {\it Every polynomial pointed out in the consequence has the
degree
$m$, as the polynomial of $X$.     }
\vskip 6mm

\noindent
\rm
3.6. \ The sequence $a_i,b_i =a_i^{-1}$, $i=1,\infty$, $a_i\in G$
of the theorem 3.5 gives the map $B_{\infty} \to B_{\infty}(G)$,
such that $\tilde t_i^2 =1$. Thus the invariants of the theorem 3.5 (and
the consequence 3.5) do not differ links $\hat \beta_1 $ and
$\hat \beta_2 $ if the permutations that define the braids $\beta_1$ and
$\beta_2$
coincide (despite this, they may give some useful
information, for instance, about the braid index of the link, especially in
the case when $G$ is a free group freely generated by the elements $a_i$).

To avoid this, we have to consider the case, when $a_ib_i\ne 1$.
\vskip 6mm

\rm
3.7. {\bf Theorem.}  \ {\it For the field $k$ and any $m\in N$ let
$K=k ( \{ a_{ij}^{(s)}\} )$ be the field of the
rational functions of the undeterminates $\{a_{ij}^{(s)}\}$,
$i,j=1,...,m$, $s\in N$ over the field $k$.
Let $K[T,T^{-1}]$ be the ring of Laurent polynomials of the undeterminate
$T$ over the field $K$. For the sequence $a_s=(a_{ij}^{(s)} )$,
$b_s=T\cdot (a_{ij}^{(s)})^{-1} $ of the
elements of $GL_m (K[T,T^{-1} ])$ and every link $\hat \beta $
($\beta \in B_n$ for some $n\in N$) consider the family of
characteristic polynomials of any integer power $t$ of
  representatives of the conjugacy
classes, indicated in the theorem 3.5. For $0\le l < m$
let us consider the product $P_l =T^{t(l-m) \exp \beta }$.
$P_1^{(l)} ... P_r^{(l)}$, where $r$ is the number of knots in link $\hat \beta
$,
$P_1^{(l)}, ... , P_r^{(l)}$ are the multipliers nearby $X^l$ of these
characteristic polynomials, $\exp \beta $
is the exponent sum of $\beta $
($\exp \beta = \alpha_1 +\cdots + \alpha_p$ if $\beta=t_{i_1}^{\alpha_1} ...
t_{i_p}^{\alpha_p})$.
Then the sequence $P_0,...,P_{m-1}$ of elements in $K[T,T^{-1}]$
is invariant of the link $\hat \beta $.   }
\vskip 3mm

\rm
{\bf Proof.} \ The fact that the sequence  $P_0,...,P_{m-1}$ does not
depend of the Markov move  3.4 (i) is obvious. If we have the
Markov move 3.4 (ii), then one and only one of conjugacy classes in
$GL_m (K[T,T^{-1}])$  that are given in theorem 3.5
multiplies on $T$ if the move is
$\beta \equiv t_n \beta $, or it multiplies on $T^{-1}$, if the move is
$\beta \equiv t_n^{-1}\beta $. Taking into account that
$\exp ( t_n \beta)=\exp \beta +1$,
$\exp (t_n^{-1}\beta ) =\exp \beta -1$ we obtain that the sequence
$P_0,...,P_{m-1}$ does not depend of Markov moves and thus
it is invariant of the link $\hat \beta $.
\vskip 6mm

{\bf Remark.} \ {\it Let $\hat \beta $ be a link. Any $m\in N$
and any integer $t$ define sequence $P_0,...,P_{m-1}$, and    }
$P_i \in k\left[ \{a_{ij}^{(s)}\}, {1\over {\rm det}\,( a_{ij}^{(s)}) }, X,
T,T^{-1}\right]$.
\vskip 6mm

\rm
3.8. {\bf Theorem.}  \ {\it Let $G$ be the group and $Tr :G\to k$ be the map
from $G$ into
the field $k$ such that $Tr(g^{-1} x g) =Tr(x)$ for any
$g,x \in G$. Suppose that there exists an element $u\in G$ such that
$Tr(ux) =\lambda_1 Tr(x)$, $Tr(u^{-1}x) =\lambda_2 Tr(x)$,
for any $x\in G$,  where $\lambda_1 ,\lambda_2 \in k$, $\lambda_1 \lambda_2\ne
0$.
Then any sequence $a_s\in G$, $s\in N$ of elements of the group $G$
defines invariant of any link $\hat \beta $ $(\beta \in B_n$ for some
$n$) obtained in the following way. For the sequence
$a_s,b_s =a_{s-1} ... a_1 u a_1^{-1} ... a_s^{-1} a_s^{-1}$
$s=1,2...$ (this sequence satisfies the relations 3.3 (i)) we consider
the representatives $\gamma_1,...,\gamma_r$ of conjugacy classes for
$\hat\beta$, that is given by the theorem 3.5. Let
\begin{displaymath}
\matrix { ({\rm i}) &   I_{\hat\beta}= k_n \cdot V^{\exp \beta }\cdot
{\prod \limits_{i=1}^{r} Tr(\gamma_i)}  \hfill\cr },
\end{displaymath}
where  $r$ is the number of knots in link $\hat\beta $,
\begin{displaymath}
V^2 =\lambda_2 \cdot \lambda_1^{-1},
\quad k_1 =1, \ k_{n+1}=k_n\cdot V^{-1} \cdot \lambda_1^{-1}.
\end{displaymath}

Then element $I_{\hat\beta} \in k\left( \sqrt{ {\lambda_2\over \lambda_1}}
\right)$
is invariant of the link $\hat\beta $.   }
\vskip 3mm

\rm
{\bf Proof.} \ It is analogous to the proof of the theorem 3.7.
The only difference is that we consider invariant in the form (i). The
condition $I_{\widehat{t_n\beta }} =I_{\widehat{t_n^{-1}\beta }}$
implies equality $\lambda_1 V {\prod \limits_{i=1}^{r} }Tr(\gamma_i)=
\lambda_2 V^{-1} {\prod \limits_{i=1}^{r} }Tr(\gamma_1)$
and hence $V^2 =\lambda_2 \lambda_1^{-1}$. The condition
$I_{\hat\beta } =I_{{\widehat{t_n\beta }}}$ implies equality
\begin{displaymath}
k_n\cdot {\prod \limits_{i=1}^{r} }Tr(\gamma_i)= k_{n+1} \cdot V\cdot
\lambda_1 \cdot {\prod \limits_{i=1}^{r} }Tr(\gamma_i)
\end{displaymath}
that implies relation $k_{n+1} =k_n\cdot V^{-1} \cdot \lambda_1^{-1}$.
\vskip 6mm

{\bf Remark.} \  {\it In fact one can easily prove that if $Tr \not\equiv 0$,
then $\lambda_2 =\lambda_1^{-1}$ and hence $V=\pm \lambda_1^{-1}$.   }
\vskip 6mm

In the followings items we give another one approach to the
construction of invariants of oriented knots and links.
\vskip 4mm

\noindent 3.9.  \ Let $\overline {\cal B}$ be an algebra over the ring
$ Z$ with generators $A, B, C, D, A_1$,\linebreak $ B_1, C_1, D_1$
and the following defining relations:
\begin{displaymath}
\matrix{
({\rm i})    &  A^2 +BAC =A,\hfill\cr
({\rm ii})   & D^2 +CDB=D,\hfill\cr
({\rm iii})  & CB-BC =ADA -DAD,\hfill\cr
({\rm iv})   & BA- AB=BAD,\hfill\cr
({\rm v})    & AC-CA =DAC,\hfill\cr
({\rm vi})   & DB-BD =ADB,\hfill }
\matrix{
({\rm vii})  & CD -DC =CDA,\hfill\cr
({\rm viii}) & AA_1 +BC_1 =1,\hfill\cr
({\rm ix })  & CA_1 +DC_1 =0,\hfill\cr
({\rm x})    & AB_1 +BD_1 =0,\hfill\cr
({\rm xi})   & CB_1 +DD_1 =1.\hfill\cr }
\end{displaymath}

Then $\overline{{\cal B}} ={\cal B}\langle A_1 ,B_1, C_1, D_1 \rangle$,
where   ${\cal B}$
 is the braid algebra (if we replace in the definition of the braid
algebra the field $k$ on the ring $ Z $).

Note that in the group $GL_2 (\overline{{\cal B}})$
matrix
$\left(\matrix{ A & B\cr C & D}\right) $  is invertible and
\begin{displaymath}
({\rm xii})\quad \left(\matrix{
A   & B   \cr C   & D   }\right)^{-1} =\left( \matrix{ A_1 & B_1 \cr C_1 & D_1
}\right).
\end{displaymath}

It is the main reason to introduce the generators $A_1, B_1, C_1, D_1$.

Let $X$ be formal "noncommutative monom" of the elements
$A,A_1,B$, \linebreak $B_1,C,C_1, D,D_1$ that is the product of some of the
elements
$A, A_1, B,B_1,\linebreak  C,C_1,D,D_1$ of the algebra $\overline{{\cal B}}$ in
some order
(for instance $X=A_1^3 \,B_1 \,A_1
\linebreak B^2\,C_1^5\,D\,C_1$). Let $l(X)$
           be the length of this monom. For $1\le i\le l(X)$
\linebreak let $\psi (i) =n_1 (i) +n_2 (i) -n_3(i)-n_4(i)$,
where $n_1(i),n_2(i), n_3(i), n_4(i)$ are the numbers of the appearences of
$B,B_1, C,C_1$, respectively in the part of this monom from the first
position to the $i^{{\rm th}}$ position (for the above example
$n_1 (4)=0, n_2(4)=1, n_3(4)=0, n_4(4)=0, n_1(12)=2,
\linebreak n_2(12)=1, n_3(12)=0, n_4(12)=6)$.
\vskip 4mm

{\bf Definition.} \ {\it We shall call the subgroup in additive group
$\overline {{\cal B}}$ ge\-nerated by $1\in \overline{{\cal B}}$ and monoms $X$
(after real multiplication), for which
$\psi (l(X))=0$ $(\psi (l(X))=0$ and $\psi (i)\ge 0$, for any $1\le i \le
l(X))$
the group $T_{nc}$ of diagonal elements (the group $T_{ncl}$ of the
diagonal elements in the last position of the diagonal).  }
\vskip 5mm

Note that in reality both of this subgroups are subrings of the
ring $\overline{{\cal B}}$.
We will denote by $J$ the right ideal of     ring $T_{ncl}$ generated by
elements
 $A-1, A_1 -1$.  Thus $J=(A-1) T_{ncl} +(A_i-1)T_{ncl}$ and
$J\subset T_{ncl} \subset T{nc}$.
For any $k\in N$ let
\begin{displaymath}
\{ a_{ij}\}_{\scriptstyle i=1,k\atop j=1,k},\quad
\{ b_{ij}\}_{\scriptstyle i=1,k\atop j=1,k}, \quad
\{ c_{ij}\}_{\scriptstyle i=1,k\atop j=1,k},\quad
\{ d_{ij}\}_{\scriptstyle i=1,k\atop j=1,k}
\end{displaymath}
be
$4k^2$ undeterminates. Let
$P ={ Z} [\{ a_{ij}\},\{ b_{ij}\},\{ c_{ij}\},\{ d_{ij}\},u]$
be the commutative algebra over ${Z}$ (that is the ring) generated
by these $4k^2$ undeterminates and undeterminate $u$, with  following defining
relation:
$u \times {\rm det} \left(\matrix{ \varphi (A )& \varphi (B)\cr \varphi (C) &
\varphi (D)}\right) =1$
plus 7$k^2$ relations that come from relations (i)-(vii) if we
rely
\begin{displaymath}
\varphi (A)=(a_{ij})_{k\times k},\ \
\varphi (B)=(b_{ij})_{k\times k},\ \
\end{displaymath}
\begin{displaymath}
\varphi (C)=(c_{ij})_{k\times k},\ \
\varphi (D)=(d_{ij})_{k\times k}
\end{displaymath}
and substitute in the relations (i)-(vii)
$\varphi (A), \varphi (B), \varphi (C), \varphi (D)$ instead of
$A, B, C, D$ respectively
(thus the algebra $P$ gives the "common point" of those $k$-dimensional
representations of the  braid algebra ${\cal B}$, that lead to the
representations of the braid group $ B_n$ for any $n\in N$
according to the theorem 2.2). Matrix
$ \left(\matrix{ \varphi (A) & \varphi (B)\cr \varphi (C) & \varphi (D)}\right)
$
over the ring $P$ is invertible . Let us define $\varphi (A_1)$, $\varphi
(B_1)$,
$\varphi (C_1)$, $\varphi (D_1)$ by the equality
$\left(\matrix{ \varphi (A_1) & \varphi (B_1)\cr \varphi (C_1) & \varphi
(D_1)}\right) =
\left(\matrix{ \varphi (A) & \varphi (B)\cr \varphi (C) & \varphi
(D)}\right)^{-1}$
(i.e. matrices $A_1, B_1, C_1, D_1$ become concrete matrices with
entries from the ring $P$).
Thus we have obvious homomorphism
$\varphi : \overline{{\cal B}} \to {\rm Mat}_{k\times k} (P)$
and map ${\rm tr}\,\varphi  : \overline{{\cal B}} \to P$
that to every element  $b\in \overline{{\cal B}}$,
 corresponds ${\rm tr}\,\varphi (b)$  (where trace is the  usual trace of the
matrix, i.e. the sum of its diagonal elements). Note that
${\rm tr}\,\varphi ( T_{nc})$, ${\rm tr}\,\varphi (T_{ncl})$,
${\rm tr}\,\varphi (J)$ are subgroups of additive group of the ring
$P$ and ${\rm tr}\,\varphi (J) \subset {\rm tr}\,\varphi ( T_{ncl})
\subset {\rm tr}\,\varphi ( T_{nc})$.
\vskip 6mm

{\bf Definition.} \ {\it We shall call the subgroup in
$L_k^{{\rm link}}(L_k^{{\rm d.link}}))$ additive group $P$
that is the image
${\rm tr}\,\varphi  T_{nc} ({\rm tr}\,\varphi (J))$
of the group $T_{nc}$ (of the group
$J$) under this map ${\rm tr}\,\varphi $ --- the $k^{{\rm th}}$ group
of links (the $k^{{\rm th}}$ group
of degenerate links) and the quotient group $L_k^{{\rm inv}}=
L_k^{{\rm link}}/L_k^{{\rm d.link}}$
--- the $k^{{\rm th}}$-group  of invariants of links.   }
\vskip 6mm

\noindent
3.10. \ For any $n\in N$ let us consider the homomorphism $\pi$ of
the braid group ${ B}_n =\langle t_1,...,t_{n-1} \rangle$
in the group $GL_{nk} (P)$, such that
\begin{displaymath}
\pi (t_i) =\left( \matrix{
I_{(i-1)k\times (i-1)k} & 0  & 0\cr\cr
0                       & \matrix{\varphi (A) &\varphi (B)\cr \varphi (C) &
\varphi (D)}  &  0\cr\cr
0                       & 0  & I_{(n-i)k\times (n-i)k}  } \right),
\end{displaymath}
According to the items  2.2, 3.9 such homomorphism exists and it
is obviously unique.
\vskip 6mm

\noindent
{ \bf Lemma .} \ {\it For any $\beta \in { B}_n$
the trace of the matrix $\pi (\beta )$
belongs to the group   $L_k^{{\rm link}}$, and
$k\times k$ matrix in the right low corner  of the matrix $\pi (\beta )$
belongs to the subring $\varphi (T_{ncl})$
of the ring ${\rm Mat}_{k\times k} (P)$ $\pi $.  }
\vskip 5mm

\noindent
{\bf Proof.} \ Let us consider in $kn \times kn$ matrix $\pi (\beta )$ $n$
square $k\times k$ submatrices each with left upper vertix laying on the
diagonal
of matrix $\pi (\beta )$ and which do not intersect. Then from
the definition of the representation $\pi$ every of this matrices belongs to
the
$\varphi (T_{nc})$ and the $k\times k$ matrix in the right low corner
of matrix $\pi (\beta )$ belongs to $\varphi (T_{ncl})$.
\vskip 4mm

For any commutative ring ${\cal A}$ and any $g, g_1 \in GL_r ({\cal A})$
${\rm tr} (g_1 g g_1^{-1} )={\rm tr}\, g $
(if in the case when ${\cal A}$ is a field it is well known fact, the common
case reduces to this one
if we consider the "abstract" matrices $\tilde g_1= (x_{ij})_{r\times r}$,
$\tilde g= (y_{ij})_{r\times r}$ over the field
$Q(\{x_{ij} \}, \{y_{ij}\}$, where  $x_{ij}, y_{ij}$
$i,j =1,...,r$ are undeterminates).
\vskip 6mm

\noindent
3.11. \ For the element $l\in L_k^{{\rm link}}$ by $l^{{\wedge}}{}$
we denote the element $f(l) \in L_k^{{\rm inv}}$,

where $f$ is the natural homomorphism
$L_k^{{\rm link}} \to L_k^{{\rm link}} / L_k^{{\rm d.link}} = L_k^{{\rm inv }}$
(see definition 3.9).
\vskip 6mm

\noindent
{ \bf Theorem.} \ {\it Let $\hat \beta $ be the oriented link
obtained of the braid $\beta \in {B}_n$. Then the
element

\noindent
$$
I_{\beta,k} =[2 {\rm tr}\, \pi (\beta ) +(\exp \beta )
 ({\rm tr}\,\varphi( D_1) -{\rm tr}\,\varphi(D))
-n ({\rm tr}\,\varphi( D_1)  + {\rm tr}\,\varphi({\rm D})) ]^{{\wedge}}{}
$$
(where $\pi$ is standart representation of ${ B}_n$ into
$GL_{nk} (P)$, see 3.10,
$\varphi$ is standart homomorphism $\overline{{\cal B}} \to {\rm Mat}_{k\times
k} (P)$
, see 3.9)
of the group $L_k^{{\rm inv}}$ is invariant of the link $\hat\beta$. }
\vskip 4mm

{\bf Proof. } \ If $\beta_1, \beta_2, \gamma \in {\rm B}_n$ and
$\beta_2 =\gamma \beta_1 \gamma^{-1}$ then
$\pi (\beta_2 )=\pi (\gamma) \pi (\beta_1)$
\linebreak $ \pi (\gamma^{-1}$ and, hence,
${\rm tr}\, \pi (\beta_2) ={\rm tr}\,\pi (\beta_1)$.
So $I_{\beta , k}$ does not depend of the Mar\-kov move 3.4 (i).
Let $\beta_1 \in {\cal B}_n$, $\beta_2 \in {\cal B}_{n+1}$ and
$\beta_2 =t_n^{\delta} \,\beta_1$,
where $\delta =1$, or $\delta =-1$.
Let us consider the diagonal elements of the matrices
$\pi (\beta_1) \in GL_{nk} (P)$ and
$\pi (\beta_2) \in GL_{(n+1)k}(P)$. As follows from
the definition of the representation $\pi$
first $(n-1)k$ elements of the diagonal of matrices $\pi (\beta_1)$
and $\pi (\beta_2)$ coincide, $k\times k$ matrix $R$
in the right low corner of matrix $\pi (\beta_1)$ belongs to the algebra
$\varphi (T_{ncl})$ (see lemma 3.10).
The corresponding $k\times k$ matrix of matrix $\pi (t_n ^{\delta} \beta_1)$
(we mean $k\times k$ matrix in matrix $\pi (t_n^{\delta} \beta_1)$
coordinates of vertices of which coincide with coordinates
 of matrix $R$ in matrix $\pi (\beta_1)$ ) equals
$\varphi (A) \cdot R$ if $\delta =1$, and
$\varphi (A_1) \cdot R$ if $\delta =-1$.
The $k\times k$ matrix in the right low corner of the matrix
$\pi (t_n^{\delta}\beta_1)$
equals $\varphi (D)$, if $\delta =1$, and $\varphi (D_1)$ if $\delta =-1$.
We have
\begin{displaymath}
I_{t_n \beta_1, k} =[2{\rm tr}\, \pi (t_n \beta_1) +
({\rm tr}\,\varphi( D_1) -{\rm tr}\, \varphi(D)) (\exp t_n \beta_1 ) -
\end{displaymath}
\begin{displaymath}
-(n+1)({\rm tr}\,\varphi (D_1) + {\rm tr}\,\varphi( D))]^{{\wedge}}{}=
\end{displaymath}

\begin{displaymath}
=[2{\rm tr}\,\pi (\beta_1) + 2{\rm tr}\, (\varphi (A-1)) R +
2{\rm tr}\, \varphi (D) +
\end{displaymath}
\begin{displaymath}
+ ({\rm tr}\, \varphi (D_1) -{\rm tr}\, \varphi (D)) (\exp \beta_1 +1)-
\end{displaymath}
\begin{displaymath}
-(n+1) ({\rm tr}\, \varphi  (D_1) +{\rm tr}\,\varphi (D))]^{{\wedge}}{} =
\end{displaymath}

\begin{displaymath}
=[2{\rm tr} \pi (\beta_1) +({\rm tr} \varphi (D_1) -{\rm tr} \varphi (D))
\exp \beta_1 -n ({\rm tr} \varphi (D_1) +{\rm tr} \varphi (D))+
\end{displaymath}

\begin{displaymath}
+{\rm tr}\,\varphi (D_1) -{\rm tr} \,\varphi (D) -{\rm tr} \,\varphi (D_1)-
\end{displaymath}

\begin{displaymath}
-{\rm tr}\,\varphi (D) +2{\rm tr}\,\varphi (D) +2{\rm tr}\,(( \varphi (A-1)) R
)]^{{\wedge}}{}=
I_{\beta_1,k}
\end{displaymath}
(we took into account that $\exp t_{n} \beta_1 =\exp \beta_1 +1$
and that $2(\varphi (A-1)) R\in \varphi (J_)$ and, hence,
$2{\rm tr}\, (( \varphi (A-1)) R)\in L_n^{{\rm d. link}} )$.
We have
\begin{displaymath}
I_{t_n^{-1}\beta_1} =[2{\rm tr}\, \pi (t_n^{-1} \beta_1 ) +({\rm tr}\,\varphi (
D_1) -{\rm tr}\,\varphi ( D))
\,(\exp \beta_1 -1) -
\end{displaymath}
\begin{displaymath}
-(n+1) ({\rm tr}\,\varphi (D_1) +{\rm tr}\,\varphi  (D) )]^{{\wedge}}{}=
\end{displaymath}

\begin{displaymath}
=[2{\rm tr} \,\pi (\beta_1) +2{\rm tr}\, ((\varphi (A_1 -1)) R) +2{\rm
tr}\,\varphi ( D_1)+
\end{displaymath}
\begin{displaymath}
+ ({\rm tr}\,\varphi ( D_1) -{\rm tr}\,\varphi ( D)) (\exp \beta_1 -1)-
(n+1) ({\rm tr}\,\varphi ( D_1) +{\rm tr}\,\varphi ( D))]^{{\wedge}}{} =
\end{displaymath}

\begin{displaymath}
=[2{\rm tr}\, \pi (\beta_1)
+ ({\rm tr}\,\varphi (  D_1) -{\rm tr}\,\varphi ( D))\exp \beta_1 -
 n({\rm tr}\,\varphi (  D_1) +{\rm tr}\, \varphi (D))-
\end{displaymath}

\begin{displaymath}
-{\rm tr}\,\varphi ( D_1) +{\rm tr}\,\varphi ( D)
-{\rm tr}\,\varphi ( D_1) -{\rm tr}\,\varphi ( D) +
\end{displaymath}
\begin{displaymath}
+2{\rm tr}\,\varphi ( D_1) +2 {\rm tr\,} ((\varphi (A_1 -1)) R)]^{{\wedge}}{}=
                  I_{\beta_1,k}
\end{displaymath}
(we took into account that $\exp t_n^{-1}\beta_1 = \exp \beta_1 -1$
and that $2\varphi (A_1 -1)R\in J$.
Thus $I_{\beta_2,k} =I_{\beta_1,k} $ and, hence, $I_{\beta,k}$ does not depend
of Markov
move 3.4~ (ii).
So, according to Markov theorem 3.4, $I_{\beta,k}$ is
invariant of the link $\hat\beta$.

There is the  hope  that group $L_k ^{{\rm inv}}$ (or its suitable
quotient group) can   be calculated at
least for small $k$. For $k=1$, it easy to calculate,
that $L_1^{{\rm inv}} ={Z}$.
\vskip 3mm

Note that similar theorem can be easily formulated and proved
for the case $(n,1)$. It is more complex to formulate the analogous
theorem for the
case $(n,\lambda)$, where $\lambda > 2$, see 1.1, 2.14.
\vskip 5mm

\noindent
3.12. \ The next theorem and its consequence gives, seems, a very mighty
instrument for the construction of the invariants of oriented knots
and links.
\vskip 4mm

{\bf Theorem.} \ {\it Let $T(i)$, $i\in N$ be the sequence of tensors that
satisfy the condition of the theorem
2.1. Let for any $i\in N$ the
following matrices equalities hold
$(\gamma (i)^{i_1}{}_{i_2}) =\lambda_1 (i) \times I$,
$(\tau (i)^{i_1}{}_{i_2} )=\lambda_2(i) \times I$,
$\lambda_1(i) \in k$, \ $\lambda_2(i)\in k, \ \lambda_1 (i)\ne 0, \lambda_2
(i)\ne 0$,
where
\begin{displaymath}
\gamma(i)^{i_1}{}_{i_2} =
{\mathop {\sum}\limits_{j=1}^{{\rm dim}V_{i+1}} \, }T(i)^{i_1j}{}_{i_2j},
\quad
\tau(i)^{i_1}{}_{i_2} ={\mathop {\sum}\limits_{j=1}^{{\rm dim}V_{i+1}} \,}
T^{-1}(i)^{i_1j}{}_{i_2j}
\end{displaymath}
(($T^{-1}(i)^{i_1i_2}{}_{j_1j_2}) $ \ is matrix inverse to the matrix
$\ (T(i)^{i_1i_2}{}_{j_1j_2})\ $.
If $\ \lambda _2 (i)\cdot$ $ \lambda_1^{-1}(i)={\rm const}$
(i.e. $\lambda_2(i) \lambda_1^{-1} (i)$
does not depend of $i$) then
for any link $\hat \beta, \beta\in B_n$ element
$\ I_{\beta} =k_n V^{{\rm exp \beta}} {\rm Tr}\, (\beta) \in
k \left(\sqrt{\lambda_2(1) \lambda_1^{-1} (1)}\right)$
\ where \, $V=\sqrt{{\lambda_2(1) \lambda_1^{-1} (1)}}, k_1=1,$
$k_{n+1} =k_n V^{-1} (\lambda_1 (n))^{-1}$ \ is invariant of the link
$\hat \beta$
($\pi$ is the representation of the group $B_n$ determined by the tensor\'{}s
sequence $T(i)$,
see theorem 2.1).  }
\vskip 4mm

{\bf Proof.} \ Let $\beta_1 \in B_n$ be obtained from $\beta \in B_n $ by
Markov
move \linebreak 3.4~ {\rm (i)}. Then $\beta_1 =\gamma \beta \gamma^{-1}$ for
some $\gamma\in B_n$ and,
hence, $I_{\beta_1} =I_{\beta}$.

Let $\beta_1 \in B_{n+1}$ be obtained from $\beta $ by Markov move 3.4 (ii).
Then $\,\beta_1 =\beta t_n^{\delta}$, \, where $\delta =\pm1$.
As \  $\pi (\beta)\in GL_k (V_1 \otimes \cdots \otimes V_n)$ we have
$$
\pi (\beta) =\sum T^{i_1 ... i_n}{}_{j_1 ... j_n}
e^{(1)}_{i_1}\otimes \cdots \otimes e^{(n)}_{i_n}\otimes
e^{(1)j_1*}\otimes\cdots \otimes e^{(n)j_n *}.
$$
Then ${\rm tr}\, \pi (\beta) ={\mathop {\sum}\limits_ {i_1,...,i_n} }
T^{i_1,...,i_n}{}_{i_1,...,i_n}$
 and \ ${\rm tr}\, \pi (\beta t_n^{\delta}) = $
\begin{displaymath}
={\mathop {\sum}\limits_{i_1,...,i_n,j} } T^{i_1,...,i_n}{}_{i_1,...,i_n}
\times
T^{\delta}(n) ^{i_nj}{}_{i_n j} +
\end{displaymath}
\begin{displaymath}
+{\mathop {\sum}\limits_ {i_1,...,i_n, j, i_n' \ne i_n} }
               T^{i_1,...,i_{n-1}i_n}{}_{i_1,...,i_{n-1} i_n'}
t^{\delta}(n) ^{i_n' j} {}_{i_n j}=
\end{displaymath}
\begin{displaymath}
={\mathop {\sum}\limits_ {i_1,...,i_n}} \left(T^{i_1,...,i_n}{}_{i_1,...,i_n}
{\mathop {\sum}\limits_ {j}} T^{\delta} (n) ^{i_nj}{}_{i_n j} \right ) +
\end{displaymath}
\begin{displaymath}
+{\mathop {\sum}\limits_ {i_1,...,i_n,i_n'\ne i_n} }
   \left( T^{i_1,...,i_{n-1} i_n}{}_{i_1,...,i_{n-1} i_n'}
  {\mathop {\sum} \limits_{j} } T^{\delta}(n) ^{i_n' j}{}_{i_n j} \right)
\end{displaymath}
(here $T^{\delta} (n)^{i_1 j_1}{}_{i_2j_2} =T(n)^{i_1j_1}{}_{i_2j_2}$
if $\delta =1$ and $T^{-1}(i)^{i_1i_2}{}_{j_1j_2}$ if $\delta =-1$).
Thus ${\rm tr}\, \pi (\beta t_n)= \lambda_1(n) {\rm Tr}\, (\beta)$ and
${\rm tr}\, \pi (\beta t_n^{-1}) =\lambda_2 (n) {\rm Tr}\, (\beta)$.
Now it is an easy  exercise  to show that
$I_{\beta t_n} =I_{\beta t_n^{-1}} =I_{\beta}$ and, hence,
$I_{\beta}$ does not depend of Markov move 3.4 (i)-(ii).
According to Markov theorem $I_{\beta}$ is invariant of the link $\hat\beta$.
\vskip 4mm

{\bf Consequence.} \ {\it Let $V$ be a finite dimensional vector space over the
field $k$,
${\rm dim}\,V=n$. Let $T\in V^{\otimes 2} \otimes (V^*)^{\otimes 2}$
be tensor that satisfies the conditions of the theorem 2.2
(i.e. it satisfies braid equation and its matrix
$(T^{i_1i_2}{}_{j_1j_2})$ is nondegenerate).
Suppose that $(\gamma^{i_1}{}_{i_2})_{n\times n} =\alpha_1 I_{n\times n}$,
$(\tau^{i_1}{}_{i_2})_{n\times n} =\alpha_2 I_{n\times n}$,
$\alpha_1, \alpha_2 \in k$,
$\alpha_1\alpha_2 \ne 0$, where
\begin{displaymath}
\gamma^{i_1}{}_{i_2} ={\mathop {\sum}\limits_ {j=1}^{n}} T^{i_1j}{}_{i_2
j},\qquad
\tau^{i_1}{}_{i_2} ={\mathop {\sum}\limits_ {j=1}^{n}} (T^{-1}) ^{i_1 j}{}_{i_2
j}
\end{displaymath}
((($T^{-1})^{i_1i_2}{}_{j_1j_2} )$ is matrix inverse to the matrix
$(T^{i_1i_2}{}_{j_1j_2}))$. Then for any link $\hat \beta$, $\beta\in B_n$
element
$I_{\beta}= k_n V^{{\rm exp}\beta} {\rm Tr}\, \pi(\beta) \in
                    k\left[ \sqrt {\alpha_2 \alpha_1^{-1}} \right]$
where
$V=\sqrt {\alpha_2 \alpha_1^{-1}}$, $k_1 =1$,
$k_{n+1} =k_n V^{-1} \alpha_1 ^{-1}$
is invariant of the link $\hat \beta$
\linebreak ($\pi $ is the representation of the group
$B_n$ determined by tensor $T$, see theo\-rem 2.2).   }
\vskip 4mm
{\bf Proof.} \ It is obvious.
\vskip 4mm

Note that we can consider commutative algebra
$k[\{ T^{i_1 i_2}{}_{j_1j_2}\}, \alpha_1^{\pm 1},
              \linebreak \alpha_2, x, u]$ over the field $k$ with
generators $\{T^{i_1i_2}{}_{j_1j_2}\}$,
$\alpha_1 ,\alpha_2, x,u $ and defining relationships, that follows from the
braid
equation (see 2.2), conditions of the consequence (for instance,
$x^2=\alpha_2\alpha_1^{-1}$, $u\times {\rm det}\,
(T^{i_1i_2}{}_{j_1j_2})=\linebreak =1$).
 Thus we obtain the common point of the tensors
 that satisfy the conditions of the consequence
and universal invariant that belongs to this algebra.
But there is no confidence that this formal algebra can be  calculated.
So a concrete examples of such tensors that satisfy the conditions of the
theorem or the consequence are very important.

In the following items we give only one of such examples and,
hence, the example of invariant of oriented knots and links.
\vskip 6mm

\noindent
3.13 \ Without any doubt (especially keeping
in mind the results of the item 2.10) there are a huge amount of the decisions
of
the equations 2.1~ (viii) and 2.2 (viii), such that matrices
of the corresponding tensors are nondegenerate (see 2.1, 2.2). Next
theorem gives only one example of such decision.
\vskip 4mm

{\bf Theorem.} \ {\it Let $a_i =(\alpha(i)^{i_1}{}_{i_2})_{n\times n}$,
$b_i =(\beta(i)^{i_3}{}_{i_4})_{n\times n}$ be the sequence of nondegenerate
square $n\times n$ matrices over the field $k$, such that
$a_i b_i =b_{i+1} a_{i+1}$, and $V$ is a vector space over the field $k$,
${\rm dim}_k V=n$. Put $V_i=V$, $T(i) =\sum T(i)^{i_1i_2}{}_{i_3i_4} e_{i_1}
\otimes e_{i_2}\otimes e^{i_3 *} \otimes e^{i_4 *}$,
for any $i\in N$ and for some basis $\{ e_s\}$,
$s=1,...,n$, where $T(i)^{i_1i_2}{}_{i_3i_4} =\beta (i) ^{i_1}{}_{i_4}
\alpha(i)^{i_2}{}_{i_3}$. Then matrices $(T(i)^{i_1i_2}{}_{i_3i_4})$ are
nondegenerate and tensors $T(i)$ satisfy
the equation 2.1 {\rm (viii)} and thus (according to the theorem 2.1)
determine for any $m\in N$ the linear representation $\pi$ of the group $B_m$.
}
\vskip 4mm

{\bf Proof.} \ Simple calculation proves it. The fact that matrices
\linebreak $(T(i)^{i_1i_2}{}_{i_3i_4})$ are nondegenerate follows from the
fact, that
$T(i) \in \linebreak {\rm End}\,V^{\otimes 2}$, is induced by the linear map
$V\times V \to V\times V$ with matrix \linebreak $\left( \matrix{ 0 & b_i\cr
a_i & 0}\right)$,
where \,$a_i$ \,and $b_i$ are nondegenerate matrices.
\vskip 4mm

\noindent
3.14. {\bf Theorem.} \ {\it For the field $F$ and any $n\in N$ let $k =F(T, \{
a(i)^{i_1}{}_{j_1} \})$
be the field of the rational functions (each depends of finite
number of undeterminates) over the field $F$ of the
undeterminates $T$, $\{ a(i)^{i_1}{}_{j_1}\}$,
$i_1,j_1 =1,...,n, \ i\in N$.
 Then the sequence $a_i = (a(i)^{i_1}{}_{j_1})_{n\times n}$,
 $b_i = T \times (a(i)^{i_1}{}_{j_1})^{-1}_{n\times n}$ square $n\times n$
matrices over the field $k$ satisfies the condition of the theorem 3.13.
The sequence of  tensors $T(i)$, constructed according to the theorem 3.13,
satisfies the condition of the theorem 3.12 with
$\lambda_1 (i) =T$, $\lambda_2(i) =T^{-1}$.
Thus the sequence of tensors $T(i)$ determines, according to the theorem 3.12,
invariant of oriented knots and links. Namely, for any link
$\hat \beta$, $\beta \in B_n$ the element
 $I_{\beta}  =(T)^{-\exp \beta } {\rm Tr }\,\pi (\beta) $
 \,of the field
$k$ is invariant of the link $\hat\beta $. }
\vskip 4mm

{\bf Proof.} \ Every statement is evident.
\vskip 4mm

The approach above is considered only for the case
case $(n,\lambda )=(2,1)$ (see 2.1) but it can be easy generalized
to the case of arbitrary $(n,\lambda )$ (especially for "small"
$n,\lambda$), $n,\lambda \in N$, $\lambda < n$ (see 1.1, 2.14).
But it refers to all this paper.
\vskip 4mm

{\bf Remark.} \ {\it In a few days before sending this paper to publication
the author revealed the work [A], where the algebra ${\cal B}_{2,1}$
appeared also. In this article this algebra is considered mostly
of physical point of view and some of its unitary infinite
dimensional representation are constructed. The theorem similar to the
theorem 2.2 in part, corresponding the algebra ${\cal B}_{2,1}$,
hovewer, is not formulated, but is indicated.

\newpage

\vspace*{20mm}

\centerline {{\large Reference } }
\vskip 5mm

[GP] {\it Gelfand I.M., Ponomarev V.A.} Irreducible representations of the
Lorentz
group. Uspechi Mat. Nauk., 140 (1968), 3--60.

[NRSB] {\it Nazarova L. A., Roiter A. V., Sergeitchuk V. V.,
Bondarenko V. M.} Application of the theory of modules
over diades to the classification of finite $p$-groups
with abelian subgroup of index $p$, and to the classification of the pairs of
annihilating operators. Notes of the seminars in LOMI, USSR, 28 (1972), 69--92.

[B] {\it Birman J.S.} Braids, links and mapping class groups. Ann. Math
Studies \#82, Princeton university Press 1975.

[D1] {\it Drozd Y.A.} Tame and wild matrix problems.
Collection of articles "Representations and
quadratic forms" (in Russian).
 Institute of math. of AN Ukr. SSR, Kiev (1979), 39--73.

[D2] {\it Drozd Y.A.} Representations  of commutative algebras.
Func. analyses and its applications, 6 (1972), 41--43.

[Dr] {\it Drinfel\'{}d V.G.} Quantum groups. Proceedings of the ICM,
Berkerley, 1986, 798--820.

[J] {\it Jones V.F.R.} Hecke algebra representations of braid groups and link
polynomials. Ann. of Math. (2), 126 (1987), no.2, 335--388.

[F--O] {\it Freyd P., Yetter D., Hoste J., Lickorish W., Millet K. and
Ocneanu A.} A new polynomial invariant of knots and links.
Bull. A.M.S. 12(1985), 239--246.

[M] {\it Morton H.R.} Problems. In Proceedings of the
Joint Summer Conference on Artin\'{}s Braid Group. Santa Cruz, California
(1986), 557--574.

[L] {\it Lehrer G.I.} A survey of Hecke algebras and the Artin braid groups.
In Proceedings of the Joint
Summer Conference on Artin\'{}s Braid Group. Santa Cruz, California
(1986), 365--385.

[T] {\it Turaev V.G.} The Yang-Baxter equations and invariantsof
links. Invent. Math. 92 (1988),  pp. 527--553.

[K] {\it Kontsevich M.} Vassiliev\'{}s knot invariants, Adv. Sov. Math.,
16, part 2, (1993),137--150.

[A] {\it Arik M., Aydin F., Hizel E., Kornfilt J., Yildiz A.}
Braid group related algebras, their representations and generalized
hydrogenlike spectra. J. Math. Phys. 35 (6), June 1994, pp. 3074--3088.
\vskip 25mm

\noindent
\begin{minipage}[t]{55mm}
Department of Mathematics\\
Kiev University\\
Vladimirskaja Str, 64\\
Kiev 252617\\
UKRAINE
\end{minipage}

\end{document}